%% file: main.tex
\documentclass[11pt]{article}
\usepackage[margin=1in]{geometry}
\usepackage{amsthm, amsopn, bm, xcolor, amsmath}
\usepackage{multirow}
\usepackage{booktabs}
\usepackage{pifont}
\usepackage{hyperref}

\input{defs}

\usepackage{microtype}
\usepackage{graphicx}
\usepackage{subfigure}
\usepackage{booktabs} 

\usepackage{amssymb}
\usepackage{mathtools}
\usepackage[utf8]{inputenc}

\renewcommand{\epsilon}{\varepsilon}
\usepackage{cleveref}

\title{Signal-Plus-Noise Decomposition of Nonlinear Spiked Random Matrix Models} 

\author{
Behrad Moniri\;
and\; Hamed Hassani\footnote{The authors are with the Department of Electrical and Systems Engineering, University of Pennsylvania. Correspondence to \texttt{bemoniri@seas.upenn.edu}.}
}
\date{}

\def\hmath$#1${\texorpdfstring{{\rmfamily\textit{#1}}}{#1}}

\usepackage{setspace}

\begin{document}
\maketitle

\begin{abstract}
    In this paper, we study a nonlinear spiked random matrix model where a nonlinear function is applied element-wise to a noise matrix perturbed by a rank-one signal. We establish a signal-plus-noise decomposition for this model and identify precise phase transitions in the structure of the signal components at critical thresholds of signal strength. To demonstrate the applicability of this decomposition, we then utilize it to study new phenomena in the problems of signed signal recovery in nonlinear models and community detection in transformed stochastic block models. Finally, we validate our results through a series of numerical simulations.
\end{abstract}

\section{Introduction}
\label{sec:introduction}

Spiked Random Matrix models, where a low-rank signal is added to a large-dimensional noise matrix, have been studied extensively in recent years. Among many applications, these models have successfully been used for a theoretical study of spectral methods for signal recovery and detection (see, e.g., \cite{johnstone2001distribution, BBP, el2018estimation, perry2018optimality, barbier2018rank, el2020fundamental, guionnet2022low, pak2024optimal}, etc.), and community detection in stochastic block models (see, e.g., \cite{abbe2018communitydetection, han2023spectral, lee2023phase, mcsherry2001spectral}, etc.). In these models, for $n\in \mathbb{N}$, it is assumed that we observe 
\begin{align}
\label{eq:linear_signal_plus_noise}
    \bY_{\rm Lin} = \frac{1}{\sqrt{n}}\bW + \lambda \vx \vx^\top,
\end{align}
where $\mathbf{W} \in\R^{n\times n}$ is a noise matrix, $\vx\in \R^n$ is the unit norm signal, and $\lambda>0$ is the signal strength. 

Given the noisy observation $\bY_{\rm Lin}$, the top eigenvector of $\bY_{\rm Lin}$ is used as an estimate of $\vx$. This algorithm is based on the fact that the spectrum of $\bY_{\rm Lin}$ undergoes a sharp phase transition when varying $\lambda$: below a critical threshold of the signal strength, there is no alignment between the signal $\vx$ and the top eigenvector, and the top eigenvector is not useful for signal recovery. However, beyond this threshold, the top eigenvector exhibits a non-trivial alignment with the signal, enabling signal recovery. See Section~\ref{sec:linear_wigner} for more details.




In this paper we study an element-wise transformation of this model, which we call a \textit{nonlinear spiked random matrix model}. We define
\begin{align}
    \label{eq:main-form}
    \bY = \frac{1}{\sqrt{n}}f\left(\bW + \lambda\sqrt{n} \vx \vx^\top\right),
\end{align}
where $f: \R \to \R$ is a nonlinear function applied element-wise. This is a natural extension of spiked random matrix models. The function $f$ can model a nonlinear sensing mechanism or a function applied to truncate, censor, or normalize the data. 

We study such random matrix models in the regime where the size $n \to \infty$, and $\Vert\bW\Vert_{\rm op} = \Theta_\sP(\sqrt{n})$ and $\Vert\vx\Vert_2 \overset{P}{\to}1$. We consider a  signal strength that can potentially grow with $n$, given by
$\lambda = c_\lambda n^\alpha$ where $c_\lambda > 0$ and $\alpha \in [0, 1/2)$ and characterize the spectral properties of this model and the problem of estimating $\vx$ given the nonlinear noisy observation $\bY$. 

\subsection{Examplar Applications}
\paragraph{Transformed Stochastic Block Models.}
\label{sec:motivations}
Assume we have a stochastic block model with \( n \in \mathbb{N} \) nodes, where each node belongs to one of two communities, labeled \( -1 \) or \( 1 \). Let the vector \( \vu \in \{-1, 1\}^n \) represent the community membership of each node. We consider a weighted, fully connected graph where the weights of edges between nodes within the same community are drawn from the distribution \( \mathcal{D} \) with mean \( \gamma \), and the weights of edges between nodes in different communities are drawn from the distribution \( \bar{\mathcal{D}} \) with mean \( \bar{\gamma} \). Due to censoring, preprocessing, or the use of a nonlinear measuring device, instead of observing the weight matrix \( \mathbf{A} \), we observe an element-wise transform of it, \( f(\mathbf{A}) \), for a function \( f: \mathbb{R} \to \mathbb{R} \).

Assuming that $\gamma +\bar\gamma = 0$, the expected value of the adjacency matrix $\bA \in \R^{n \times n}$ is $\E \bA = \frac{\gamma - \bar\gamma}{2}\vu\vu^\top$, and the observed transformed adjacency matrix can be written as
\begin{align*}
    \frac{1}{\sqrt{n}}f(\bA)  = \frac{1}{\sqrt{n}} f\left((\bA - \E \bA) + \frac{\gamma - \bar\gamma}{2}\vu\vu^\top\right),
\end{align*}
which is in the form of  nonlinear spiked random matrix models of \eqref{eq:main-form} with $\bW = \bA - \E\bA$, $\vx = \vu/\sqrt{n}$, and $\lambda = \sqrt{n}(\gamma - \bar\gamma)/2$. In this problem, the goal is to recover the community membership information in $\vu$ from the observation $f(\bA)/\sqrt{n}$. We will study this problem in detail in Section~\ref{sec:pre-transformed-block}.

\paragraph{Signed Signal Recovery.} Recovering the signal $\vx$ given the nonlinear observation of \eqref{eq:main-form} can be significantly more challenging than recovering $\vx$ given the linear observation in \eqref{eq:linear_signal_plus_noise}. For instance, consider a setting where $f$ is an even function. Due to the symmetries in the problem, it becomes impossible to recover the sign of the vector $\vx$, regardless of how large the signal strength is. For a general function $f$, we show that there can be a gap in the signal strength required for signed and unsigned signal recovery and signed signal recovery can require an asymptotically larger signal strength. We will study this problem in detail in Section~\ref{sec:signed-signal-recovery}.

\paragraph{Other Examples.}  Beyond the problems discussed above, the model in \eqref{eq:main-form} appears naturally in other contexts as well. For example, \cite{ba2022high, moniri_atheory2023, cui2024asymptotics} study a two-layer neural network and demonstrate that the first layer weights, after the initial training step, can be written as a spiked random matrix and the feature matrix of the two-layer neural network conforms to the form in \eqref{eq:main-form}, with the step-size acting as the signal strength. Also, the theoretical model for unsigned low-rank matrix completion analyzed in \cite{unsigned-matrix-recovery} is a special case of such nonlinear random matrix models with $f(x) = |x|$.

\subsection{Contributions}

\begin{itemize}
    \item In Section~\ref{sec:signal+noise}, we derive a novel signal-plus-noise decomposition for $\bY$ and show that $\bY$ can be approximated in operator norm by $\frac{1}{\sqrt{n}}f(\bW)$ (i.e., noise) plus several terms depending on $\vx$ (i.e., signal). The number and structure of the signal components in this expansion depend on the power $\alpha$ and undergo a series of sharp phase transitions as $\alpha$ varies. Specifically, we show that as $\alpha$ increases from $0$ to $1/2$, new higher degree signal components emerge at critical values $\alpha_\ell = \frac{\ell - 1}{2\ell}$, $\ell \in \mathbb{N}$.
    
    \item In Section~\ref{sec:signed-signal-recovery}, we use this decomposition to study the signed signal recovery problem and characterize the precise signal strength required to recover the signal $\vx$ given $\bY$. We show that recovering the sign can require an asymptotically larger signal strength compared to the signal strength required for recovering $\vx$ up to ambiguity in its sign.  We thoroughly analyze this gap and study under what conditions such a gap exists.

    \item In Section~\ref{sec:pre-transformed-block}, we study the transformed stochastic block model problem. We fully characterize the phase transitions for weak recovery of community membership of nodes as a function of $\gamma - \bar\gamma$. Specifically, we examine how the signal strength required for recovering the community structure depends on the function $f$ and the distributions $\cD$ and $\bar\cD$.
    
    \item In Section~\ref{sec:experiments}, we confirm our theoretical findings through various numerical simulations.
\end{itemize}


\subsection{Prior Works
on Nonlinear Spiked Random Matrix Models}
Spectrum of random matrices in the form of \eqref{eq:main-form} have been studied extensively when $\lambda = 0$ (see e.g., \cite{el2010spectrum,el2010information,couillet2016kernel,pennington2017nonlinear,peche,wang2022spectral}, etc.) and has found many applications in the analysis of the performance of two-layer neural network with first layer fixed at initialization, and the second layer trained using ridge regression, i.e., random features models (see e.g., \cite{louart2018random,goldt2022gaussian, mei2022generalization,lin2021causes,adlam2019random,adlam2020neural,tripuraneni2021covariate, disagreement,loureiro2021learning}, etc.).

For the linear model of \eqref{eq:linear_signal_plus_noise}, when $\bW$ has i.i.d. non-Gaussian entries (up to symmetric constraints), \cite{lesieur2015mmse, perry2018optimality} show that there exists a nonlinear element-wise transformation of $\bY_{\rm Lin}$ that exploits this non-Gaussianity and can improve the performance of signal recovery and detection using the top eigenvector. However, they assume i.i.d. entries,  only study the problem in the setting where $\lambda = \Theta(1)$ (i.e., $\alpha = 0$), and only for a specific nonlinear function $f$.

The case where $\bW$ is a symmetric matrix with i.i.d. Gaussian entries has been studied for a general function $f$ in \cite{guionnet2023spectral,feldman2023spectral} in a signal strength regime where the spectrum of eigenvalues of $\bY$ consists of a bulk of eigenvalues that stick together and only one outlier (spike). In the deep learning theory literature, and in the context of feature learning, \cite{moniri_atheory2023} studies nonlinear spiked random matrix models in the case where $\bW$ that is the product of a symmetric random matrix with i.i.d. Gaussian entries and an independent matrix with rows drawn uniformly from the unit sphere.



In applications such as transformed stochastic block models, the entries of the noise matrix are not identically distributed and the signal-plus-noise decompositions from prior work do not hold. Also, for signed signal recovery, the natural scaling of the signal strength should lead to several outlying eigenvalues in the spectrum and again it cannot be handled by prior work.

In this paper, we provide a general purpose signal-plus-noise decomposition that does not assume Gaussianity, or identical distribution of the noise matrix entries. The signal-plus-noise decomposition of Theorem~\ref{thm:sig+noise} has a novel structure that  shows how the distribution of the entries of $\bW$ and the nonlinear function $f$ interact and affect the spectrum of $\bY$. Also, we provide our results for a signal strength $\lambda = c_\lambda n^\alpha$, for all $\alpha \in [0,1/2)$ and $c_\lambda \in \R$. This enables us to study problems in regimes where there are multiple outlying eigenvalues which is essential for the analysis of signed signal recovery.

\subsection{Notations}
\label{sec:notations}
For $d \in \mathbb{N}$, we denote $[d] = \{1,\ldots,d\}$.
Throughout the paper, $\Theta(\cdot)$, $\Omega(\cdot)$, $\omega(\cdot)$, $O(\cdot)$ and $o(\cdot)$ correspond to standard asymptotic notations and use subscript $\sP$ for the same notions holding in probability. The symbol $\overset{\mathrm{P}}{\to}$ denotes convergence in probability. Given a function $f:\R\to\R$, for any $k \in \mathbb{N}$, we denote $f^{(k)}$ to denote the $k$-th derivative of $f$ and let $f^{(0)} = f$. For two matrices $\mathbf{A}$ and $\mathbf{B}$ of the same size, $\bA\circ \bB$ denotes their Hadamard product and for a non-negative integer $k$, $\mathbf{A}^{\circ k}$ is the matrix of the $k$-th powers of the elements of $\mathbf{A}$. We denote the standard deviation of a random variable $f(Z)$ with $Z\sim \cD$ with $\mathrm{SD}_{Z\sim \cD}(Z)$.

\section{Review of Spiked Wigner Models}
\label{sec:linear_wigner}
In this section, we will review the classical results for the (linear)
spiked Wigner model. First, recall the definition of Wigner random matrices.
\begin{definition}
    \label{def:Wigne/Wigner-Type}
    Let $\bW \in \R^{n \times n}$ be a symmetric random matrix. We call $\bW$ a Wigner matrix and write $\bW \sim \mathrm{Wig}_n(\cD)$ if its elements  are independent up to symmetric constraints and are identically distributed according to a distribution $\cD$ which is assumed to have bounded fifth moment. 
\end{definition}

Let $\cD$ be a distribution on $\R$ with mean zero and variance $\sigma_w^2$. In the spiked random matrix of \eqref{eq:linear_signal_plus_noise}, let $\bW \sim \mathrm{Wig}_n(\cD)$, and assume that $\Vert\vx\Vert_2 \overset{\mathrm{P}}{\to} 1$. It is known that with probability one, the empirical eigenvalue distribution of $\bW/\sqrt{n}$ will converge to the semi-circle law which is supported on $[-2\sigma_w, 2\sigma_w]$. See \cite[Theorem 2.5]{bai2010spectral}. The spectrum of $\bY_{\rm Lin}$ has been studied extensively in the literature (see e.g., \cite{peche2006largest, benaych2011fluctuations,benaych2011eigenvalues,capitaine2012central,bloemendal2013limits,pizzo2013finite}, etc.). It is known that when $\lambda < \sigma_w$, the top eigenvalue is again $2\sigma_w$ and its corresponding eigenvector will not be correlated with $\vx$. If $\lambda>\sigma_w$, the top eigenvalue will be isolated and larger than $2\sigma_w$ and the corresponding eigenvector to the top eigenvalue will have a non-trivial alignment with $\vx$. The following classic theorem, often called the BBP phase transition, formalizes this phase transition of the leading eigenvalue of $\bY_{\rm Lin}$. See e.g., \cite[Theorem 1.1]{pizzo2013finite}.
\begin{theorem}[BBP Phase Transition]
    \label{thm:BBP}
    Let $\lambda_1, \vu_1$ be the leading eigenvalue and its corresponding eigenvector of the matrix $\bY_{\rm Lin}$. Assuming that $\Vert\vx\Vert_2 \overset{\mathrm{P}}{\to} 1$, if $\lambda < \sigma_w$, we have:
    \begin{itemize}
        \item If $\lambda_1 \overset{\mathrm{P}}{\to} 2 \sigma_w$, and $\langle\vu_1, \vx\rangle^2 \overset{\mathrm{P}}{\to} 0$.
        \item If $\lambda > \sigma_w$, we have $\lambda_1 \overset{\mathrm{P}}{\to} \lambda + \frac{\sigma_w^2}{\lambda}$, and $\langle\vu_1, \vx\rangle^2 \overset{\mathrm{P}}{\to} 1 - \frac{\sigma_w^2}{\lambda^2}$.
    \end{itemize}
    
\end{theorem}




\section{Signal-Plus-Noise Decomposition}
\label{sec:signal+noise}
In this section, we will provide a {signal plus noise} decomposition for the nonlinear spiked random matrix model in different regimes of the signal strength $\lambda$. This decomposition will be the key tool that we will use to characterize the spectrum of $\bY$ and to analyze the performance of spectral algorithms applied to it.

Our theoretical analysis applies under the following conditions.
\begin{condition}
\label{cond:strength} We make the following assumptions on the signal $\vx$ and the noise matrix $\bW$.
    \begin{itemize}
        \item The matrix $\bW$ satisfies $\Vert\bW\Vert_{\rm op} = \Theta_\sP(\sqrt{n})$ and its elements have bounded moments.
        
        \item The $\ell_2$-norm of $\vx$ satisfies $\Vert\vx\Vert_{\rm 2} \overset{\mathrm{P}}{\to} 1$ and the $\ell_\infty$-norm of $\vx$ satisfies $\Vert\vx\Vert_{\infty} = O_\sP\left(\kappa(n)/\sqrt{n}\right)$, where $\kappa(n)$ is $O_\sP(\mathrm{polylog}(n))$. 
    \item We assume that for any $k \in \mathbb{N}$, we have $\left\Vert\vx^{\circ k}\right\Vert_{2} = \Theta_\sP({n^{(1-k)/2}})$.
    
    \end{itemize}
\end{condition}
This condition is fairly general and properly normalized subgaussian vectors $\vx$ and properly normalized matrices $\bW$ with independent subgaussian entries satisfy these conditions (see e.g., \cite[Chapter 3 and Chapter 4]{vershynin2018high}). The conditions on $\vx$ hold also when the entries of $\vx$ are i.i.d. sub-Weibull random variables (see e.g., \cite{vladimirova2020sub,sambale2023some}, etc.). We also make the following assumption of the function $f$. 

\begin{condition}[Polynomial Bound]
\label{cond:poly-bound}
    We assume that the function $f$ is infinitely differentiable almost everywhere, and there exists $C_f\in \R$ and $m_f \in \mathbb{N}$ such that
    \begin{align*}
        |f^{(k)}(x)| \leq C_f + x^{m_f}, \quad \text{for all }\quad k \in \mathbb{N},
    \end{align*}
    whenever $f^{(k)}(x)$ exists.
\end{condition}
Note that, for example, any polynomial with a finite degree satisfies this condition. Under these conditions, we can state the main theorem of the paper.
\begin{theorem}
    \label{thm:sig+noise}
    Let $\lambda = c_\lambda n^\alpha$ with $\alpha \in [\frac{\ell-1}{2\ell},\frac{\ell}{2\ell+2})$ for some $\ell \in \mathbb{N}$, and $c_\lambda \in \R$. If Conditions \ref{cond:strength}--\ref{cond:poly-bound} hold, we have $\Vert \bY - \tilde \bY\Vert_{\rm op} = o_\sP(1)$, where
    \begin{align*}
        \tilde\bY =  \frac{1}{\sqrt{n}}f(\bW) + \sum_{k = 1} ^{\ell} \bH_k,\quad \text{with}\quad\bH_k = \frac{\lambda^k n^{(k-1)/2}}{k!}\left( \E f^{(k)}(\bW)\right) \circ \left(\vx^{\circ k}\vx^{\circ k\top}\right),
    \end{align*}
where the function $f$ and its derivatives are applied element-wise.
\end{theorem}
The proof of this theorem can be found in Section~\ref{sec:thm:approximation-proof}. This theorem shows that if the signal strength grows as $\lambda = c_\lambda n^\alpha$ with $\alpha \in [\frac{\ell-1}{2\ell},\frac{\ell}{2\ell+2})$, the matrix $\bY$ can be approximated in operator norm by the matrix $f(\bW)/\sqrt{n}$ plus the sum of $\ell$-many terms that contain information about the signal $\vx$. The rank of each of these terms $\bH_k$ can be upper bounded as $\mathrm{rank}(\bH_k) \leq \mathrm{rank}(\E f^{(i)}(\bW))$.

Working directly with the nonlinear model in \eqref{eq:main-form} is difficult because of the existence of interactions between the signal $\vx$ and the noise $\bW$ through the nonlinearity $f$. Theorem~\ref{thm:sig+noise} shows that in fact the model is equivalent to a spiked matrix where the noise and signal terms are disentangled. To demonstrate the applicability of this decomposition, in the next section we will apply it to two natural nonlinear signal recovery problems.

\section{Applications}

\subsection{Signed Signal Recovery}
\label{sec:signed-signal-recovery}
In this section, we consider a general function $f$ and use the decomposition of Theorem~\ref{thm:sig+noise} to study the problem of \textit{signed} signal recovery described in Section~\ref{sec:motivations}. 

We first derive a signal-plus-noise decomposition of the nonlinear spiked matrix model in \eqref{eq:main-form}.
\begin{proposition}
    \label{prop:wigner}
    Let $\bW \sim \mathrm{Wig}(\cD)$ and denote $\mu_{f^{(k)}} = \E_{Z \sim \mathcal{D}} f^{(k)}(Z)$, for any $k \in \mathbb{N} \cup \{0\}$. Let $\lambda = c_\lambda n^\alpha$ with $\alpha \in [\frac{\ell-1}{2\ell},\frac{\ell}{2\ell+2})$ for some $\ell \in \mathbb{N}$, and $c_\lambda \in \R$. If Conditions \ref{cond:strength}--\ref{cond:poly-bound} hold, we have
    \begin{align*}
        \bY = \mathbf{\Lambda} +  \frac{\mu_{f^{(0)}}}{\sqrt{n}} \mathbf{1}\mathbf{1}^\top +  \sum_{k = 1} ^{\ell} \frac{\lambda^k n^{(k-1)/2}\mu_{f^{(k)}}}{k!}\,\vx^{\circ k}\vx^{\circ k\top} + \bDelta,
    \end{align*}
where $\mathbf{\Lambda} = \frac{1}{\sqrt{n}}\left(f(\bW) - \E f(\bW)\right)$  is a mean-zero Wigner matrix, and $\Vert\bDelta\Vert_{\rm op} = o_\sP(1)$.
\end{proposition}
The proof for this theorem can be found in Section~\ref{sec:thm:wigner}. This proposition states that when $\lambda = c_\lambda n^\alpha$ with $\alpha \in [\frac{\ell-1}{2\ell},\frac{\ell}{2\ell+2})$, the spectrum of the matrix $\bY$ will consist of the bulk spectrum resulting from $\mathbf{\Lambda}$, and up to $\ell + 1$ rank-one spikes.  These spikes are aligned to the element-wise powers of the vector $\vx$ and the number of these spikes is non-decreasing in $\alpha$. 
This signal-plus-noise decomposition shows that a spectral algorithm can be used for signal recovery in this problem. Using the assumption on the $\ell_2$ norm of $\vx^{\circ k}$ in Condition~\ref{cond:strength}, for the signal strength of each signal component, we have
\begin{align}
    \label{}
    \left\Vert\frac{\lambda^k n^{(k-1)/2}\mu_{f^{(k)}}}{k!}\,\vx^{\circ k}\vx^{\circ k\top}\right\Vert_{\rm op} = \Theta_\sP\left(n^{k(\alpha - 1/2) + 1/2}\right).
\end{align}
Recall that we are assuming $\alpha<1/2$, which mean that for any fixed $\alpha$, terms with higher degree $k$ have smaller signal strength.

If at least one odd power of the vector $\vx$ appears in the expansion of $\bY$, the top-$\ell$ eigenspace of the $\bY$ will contain information about $\vx$ and a spectral algorithm can be used to recover it. However, if the rank-one terms in the expansion of $\bY$ only consist of even powers of $
\vx$, the top-$\ell$ eigenspace of $\bY$ will lose all information on the sign of $\vx$.

To focus solely on the problem of sign recovery in nonlinear matrix models, we make the following assumption on the vector $\vx$.
\begin{condition} 
\label{cond:rademacher}
We assume that the signal vector $\vx \in \R^n$ satisfies $\vx =  \boldsymbol{\zeta}/\sqrt{n} \quad \text{ with }\quad \boldsymbol{\zeta} = [\zeta_1, \dots, \zeta_n]^\top \in \R^n, \quad \text{ where }\quad \zeta_i \overset{\mathrm{i.i.d.}}{\sim} \mathrm{Rademacher}(1/2).$
\end{condition}

This vector satisfies the conditions in Proposition~\ref{prop:wigner}. Note that, the vector $\boldsymbol{\zeta}^{\circ k}$ equals $\boldsymbol{\zeta}$ if $k$ is odd and equals to $\mathbf{1} \in \R^n$ if $k$ is even. Thus, using Proposition~\ref{prop:wigner}, assuming that $\lambda = c_\lambda n^\alpha$ with $\alpha \in [\frac{\ell-1}{2\ell},\frac{\ell}{2\ell+2})$ for some $\ell \in \mathbb{N}$, we can write
\begin{align}
    \label{eq:Wigner-Rad-decomposition}
    \bY = \mathbf{\Lambda} +    \left(\sum_{k = 1, k \text{ odd}} ^{\ell} \frac{c_\lambda^k n^{(\alpha-1/2)k}\mu_{f^{(k)}}}{k! \sqrt{n}}\right)\boldsymbol{\zeta}\boldsymbol{\zeta}^{\top} + \left(\sum_{k = 0, k \text{ even}} ^{\ell} \frac{c_\lambda^k n^{(\alpha-1/2)k}\mu_{f^{(k)}}}{ k!\sqrt{n}}\right)\mathbf{1}\mathbf{1}^{\top} + \bDelta,
\end{align}
with $\Vert\mathbf{\Delta}\Vert_{\rm op} = o_\sP(1)$. To analyze the performance of spectral algorithms for signal recovery, we define the even-index and the odd-index of a function $f$ and distribution $\cD$ as follows.
\begin{definition}
    Given a distribution $\cD$ supported on $\R$ and a function $f:\R \to \R$, we define the even-index and the odd-index, denoted by $I_e$ and $I_o$ respectively, as 
    \begin{align*}
        &I_e(f;\cD):= \min \{k \in \mathbb{N} \cup \{0\} \text{ and } k \text{ even}\;|\; \mu_{f^{(k)}} \neq 0\},\,
        \text{  and  }\\ &I_o(f;\cD):=  \min \{k \in \mathbb{N} \text{ and } k \text{ odd}\;|\; \mu_{f^{(k)}} \neq 0\},
    \end{align*}
    where $\mu_{f^{(k)}}=\E_{Z \sim \mathcal{D}} f^{(k)}(Z)$. We will drop the argument $(f;\cD)$ when it is clear from context. 

\end{definition}

In the case where $I_o > I_e$, the top eigenvalue will corresponds to the eigenvector $\mathbf{1} \in \R^n$. If the signal strength is large enough, the second top eigenvector of $\bY$ will be aligned to $\boldsymbol{\zeta}$ and can be used to recover it. This is formalized in the following theorem.
\begin{theorem}
    \label{thm:IO>IE} Let $\bW \sim \mathrm{Wig}(\cD)$ and $\lambda = c_\lambda n^\alpha$ with $\alpha, c_\lambda \in \R$. Assume that Conditions~\ref{cond:strength},\ref{cond:poly-bound} and \ref{cond:rademacher} hold and $I_o > I_e$. Define
    \begin{align*}
        \kappa = \frac{c_\lambda^{I_o}}{I_o!}\left(\E_{Z\sim \cD} f^{(I_o)}(Z)\right)\quad \text{ and } \quad \sigma_f = \mathrm{SD}_{Z\sim \cD}(f(Z)).
    \end{align*}
    Then, the second top eigenvalue and eigenvector pair of $\bY$, $(\gamma_2, \vu_2)$, satisfies the following. 
    \begin{enumerate}
        \item [(a)] If $\alpha < \frac{I_o - 1}{2 I_o}$, we have  $\gamma_2 \overset{\mathrm{P}}{\to} 2\sigma_f$ and $\frac{|\langle \vu_2, \boldsymbol{\zeta}\rangle|}{\Vert\vu_2\Vert_2 \Vert\boldsymbol{\zeta}\Vert_2} \overset{\mathrm{P}}{\to} 0$.
        \item [(b)] If $\alpha = \frac{I_o - 1}{2 I_o}$, we have
        \vspace{-0.5cm}
        
        \begin{align*}
            \gamma_2 \overset{\mathrm{P}}{\to}  \begin{cases}
                2\sigma_f & \kappa < \sigma_f\\
                \kappa + \frac{\sigma_f^2}{\kappa} & \kappa > \sigma_f
            \end{cases}
            \quad\text{ and }\quad \frac{|\langle \vu_2, \boldsymbol{\zeta}\rangle|}{\Vert\vu_2\Vert_2 \Vert\boldsymbol{\zeta}\Vert_2} \overset{\mathrm{P}}{\to} \begin{cases}
                0 & \kappa < \sigma_f\\
                \sqrt{1 - \frac{\sigma_f^2}{\kappa_2^2}} & \kappa > \sigma_f
            \end{cases}
        \end{align*}
        \vspace{-0.5cm}
        \item [(c)] If $\alpha > \frac{I_o - 1}{2 I_o}$, we have  $\gamma_2 = \omega_\sP(1)$ and $\frac{|\langle \vu_2, \boldsymbol{\zeta}\rangle|}{\Vert\vu_2\Vert_2 \Vert\boldsymbol{\zeta}\Vert_2} \overset{\mathrm{P}}{\to} 1$.
    \end{enumerate} 

\end{theorem}

When $I_e < I_o$, if the signal strength is large enough, the top eigenvector will be aligned to $\boldsymbol{\zeta}$. 
\begin{theorem}
    \label{thm:IO<IE} Let conditions of Theorem~\ref{thm:IO>IE} hold, but with  $I_o < I_e$. Let $(\gamma_1, \vu_1)$ be the top eigenvalue and eingenvector pair of $\bY$. Then, Part (a), (b), and (c) of Theorem\ref{thm:IO>IE} will hold for $\gamma_1, \vu_1$ (instead of $\gamma_2$, $\vu_2)$.
    
\end{theorem}

\begin{figure}
    \centering
    \includegraphics{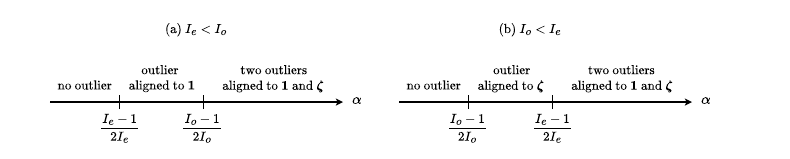}
    \caption{Emergence of outlying eigenvalues of $\bY$ with eigenvectors aligned to $\mathbf{1}$ and the signal vector $\boldsymbol{\zeta}$ in different regimes of the signal strength exponent $\alpha$ when (a) $I_e < I_o$ and (b) $I_o < I_e$.}
    \label{fig:signed}
\end{figure}

The proofs of these two theorems can be found in Section~\ref{sec:thm:IO>IE} and Section~\ref{sec:thm:IO<IE}.
These two theorems fully characterize the spectrum of $\bY$ and the performance of the spectral algorithm for the recovery of $\boldsymbol{\zeta}$ with signal strength $\lambda = c_\lambda n^\alpha$ for all choices of $c_\lambda \in \R$ and $\alpha \in [0, 1/2)$. Figure~\ref{fig:signed} summarizes these results by showing the critical values of $\alpha$ in which outlying eigenvalues with eigenvectors aligned to the vectors $\mathbf{1}$ and $\boldsymbol{\zeta}$ emerge in the spectrum of $\bY$. 

In particular, these theorems show that if for a function $f$ and a distribution $\cD$, we have $I_O > 1$, spectral algorithms asymptotically require a signal strength $\lambda = \Theta( n^{{(I_o -1)}/{2I_o}}) =  \omega(1)$ to recover the signal $\vx$ with its sign, which is asymptotically larger than the $\Theta(1)$ signal strength required in the linear case. Also, if $I_E < I_O$, the signal strength required for signed signal recovery is asymptotically larger than the signal strength required to recover the signal up to the sign. In the case where $f$ is even, we have $I_o(f;\cD) = \infty$ and signed recovery is impossible for any signal strength.

\subsection{Transformed Stochastic Block Model}
\label{sec:pre-transformed-block}
In this section, we will use the signal-plus-noise decomposition 
of Theorem~\ref{thm:sig+noise} for community detection in transformed stochastic block models described in Section~\ref{sec:motivations}.

Let $n \in \mathbb{N}$ and consider an undirected graph $\cG_n = ([n], \cE_n)$ with $\cE_n = \{(i,j)\;|\;i, j \in [n], i \neq j\}$. Assume that the vertices in this graph exhibit a community structure and they can be grouped into two non-overlapping communities $\cV_+ = \{1,\cdots,\beta n\}$ and $\cV_- = \{\beta n+1,\cdots,n\}$ for some $0<\beta<1$. Assume that the edges of this graph are weighted and define $A_{ij}$ as the weight of the edge between distinct vertices $i,j \in [n]$ that we further assume are drawn independently according to
\begin{align*}
    A_{ij} = A_{ji} \overset{\mathrm{i.i.d.}}{\sim} \begin{cases}
        \cD & i,j \in \cV_+\quad \text{or}\quad i,j \in \cV_-\\
        \bar\cD & \text{otherwise}
    \end{cases}
\end{align*}
where $\bar\cD$ and $\cD$ are two distributions over $\R$. We define $\bA \in \R^{n\times n}$ to be the adjacency matrix with entries $[\bA]_{i,j} = A_{ij}$ for all $i, j \in [n]$.  Variants of the problem of weighted stochastic block have been studied in \cite{jog2015information, li2018two, XuJogLoh}. Recently, \cite{lee2023phase} studied the problem of weak recovery in a two community weighted stochastic block model. See also \cite{jung2023detection, han2023spectral}. The problem of exact recovery in Gaussian weighted stochastic block models was also studied in \cite{pandey2024exact}.

Here we consider a transformed version of this problem where we only have access to $f(\bA)$ where $f: \R \to \R$ is a known function applied element-wise to $\bA$. The goal is to recover the community membership $\cV_0$ and $\cV_1$ using the observation $f(\bA)$.  Our results on transformed stochastic block models depend on the following quantities:
\begin{align}
    \label{def:gamma}
    &\gamma = \E_{Z\sim \cD} Z, \quad \bar\gamma = \E_{Z\sim \bar\cD} Z, \quad \gamma_{f^{(k)}} = \E_{Z \sim \cD} f^{(k)}(Z - \gamma),\quad \bar\gamma_{f^{(k)}} = \E_{Z \sim \bar\cD} f^{(k)}(Z - \bar\gamma ),\nonumber\\
    &\sigma = \mathrm{SD}_{Z \sim \cD} (f(Z-\gamma)),\quad\text{and}\quad \bar\sigma = \mathrm{SD}_{Z \sim \bar\cD} (f(Z-\bar\gamma)).
\end{align}

For simplicity, we assume that $\gamma + \bar\gamma = 0$. The expected value of  $\bA$ can be written as
\begin{align*}
    \E[\bA] = \begin{bmatrix}
        \gamma\,\mathbf{1}_{\beta n  \times \beta n} & \bar\gamma\,\mathbf{1}_{\beta n  \times (1-\beta)n}\\
        \bar\gamma\,\mathbf{1}_{(1-\beta)n \times \beta n} & \gamma\,\mathbf{1}_{(1-\beta)n \times (1-\beta)n}
    \end{bmatrix}  = \frac{\gamma - \bar\gamma}{2} \vu\vu^\top,\, \text{in which} \; \vu = 
    \begin{bmatrix}
    \mathbf{1}_{\beta n}\\[0.1cm] -\mathbf{1}_{(1-\beta)n}
\end{bmatrix}.
\end{align*}
 
Thus, $f(\bA) = f(\bar\bA + \frac{\gamma - \bar\gamma}{2} \vu \vu^\top)$ where $\bar\bA = \bA - \E \bA$ which means that $f(\bA)$ is a nonlinear spiked random matrix model defined in \eqref{eq:main-form} with $\bW = \bar\bA$, $\vx = \vu/\sqrt{n}$, and $\lambda = \frac{{\gamma - \bar\gamma}}{2} \sqrt{n}$. As a result, based on Theorem~\ref{thm:sig+noise}, we have the following proposition.
\begin{proposition}
\label{prop:blockwigner}
Let ${{\gamma - \bar\gamma}} = 2 c_\lambda n^{\alpha-1/2}$ with $\alpha \in [\frac{\ell-1}{2\ell},\frac{\ell}{2\ell+2})$ for some $\ell \in \mathbb{N}$, and $c_\lambda \in \R$. If $\bW := \bA - \E\bA$ satisfies Condition~\ref{cond:strength} and the function $f$ satisfies Condition~\ref{cond:poly-bound}, we have $\Vert\bar\bY - \frac{1}{\sqrt{n}}f(\bA)\Vert_{\rm op} = o_\sP(1)$ with
\begin{align*}
     \bar\bY = \mathbf{\Gamma} + \sum_{k = 0} ^{\ell} \frac{c_\lambda^k n^{(\alpha-1/2)k}}{k! \sqrt{n}} \left(\frac{ \gamma_{f^{(k)}} + \bar\gamma_{f^{(k)}}}{2}\vu^{\circ k}\vu^{\circ k\top}  + \frac{ \gamma_{f^{(k)}} - \bar\gamma_{f^{(k)}}}{2} \vu^{\circ (k+1)}\vu^{\circ (k+1)\top}\right),
\end{align*}
where $\mathbf{\Gamma} = \frac{1}{\sqrt{n}}(f(\bar\bA) - \E f(\bar\bA))$.
\end{proposition}
The proof can be found in Section~\ref{sec:prop:blockwigner}. 
To study the phase transitions in the eigenspace of $\bY$, we first need to study the unperturbed component $\mathbf{\Gamma}$.  In the following proposition, we apply results on Wigner-type random matrices to  characterizes the spectrum of $\mathbf{\Gamma}$ in the limit where $n\to \infty$. See Section~\ref{sec:wigner-type} for an overview of the results on the spectrum of Wigner-type random matrices. 

\begin{proposition}
\label{prop:wigner-type-block}
        Assume that $f(Z - \gamma)$ and $f(Z' - \bar\gamma)$ with $Z\sim \cD$ and $Z'\sim \bar\cD$ have bounded moments. The Stieltjes transform $m(z)$ of the empirical eigenvalue distribution of $\mathbf{\Gamma}$ satisfies 
    \begin{align*}
        m(z) = \frac{1}{n} \mathrm{tr}\left(z \bI - \mathbf{\Gamma} \right)^{-1} \overset{\mathrm{P}}{\to} \beta \bar m_1(z) + (1-\beta) \bar m_2(z),
    \end{align*}
    in which $\bar{m}_1(z)$ and $\bar{m}_2(z)$ are the solutions to the following quadratic vector equation 
    \begin{align*}
    \begin{cases}
        z \bar m_1(z) + \beta \sigma^2 \bar m_1^2(z) + (1-\beta) \bar\sigma^2 \bar m_1(z) \bar m_2(z) + 1 = 0,\\
        z \bar m_2(z) + \beta \bar\sigma^2 \bar m_1(z) \bar m_2(z) + (1-\beta) \sigma^2 \bar m_2^2(z) + 1 = 0.
    \end{cases}
    \end{align*}
\end{proposition}
The proof can be found in Section~\ref{sec:prop:wigner-type-block}. This proposition shows that the limiting spectral distribution of $\mathbf{\Gamma}$ does not carry information about the community structure. For instance, if $\beta = 0.5$, we will have $\bar m_1(z) = \bar m_2(z)$ equal to the Stiltjes transform of a semi-circle distribution.

In order for us to be able to recover the signal $\vu$ using spectral methods, an odd power of $\vu$ must appear in the expansion in Proposition~\ref{prop:blockwigner}. To characterize this, we define the signal and the constant index of a function $f$ and two distributions $\cD$ and $\bar\cD$ as follows.
\begin{definition}
\label{def:signal_index}
    Given two distribution $\cD$ and $\bar\cD$ supported on $\R$ and a function $f:\R \to \R$, we define the signal index and constant index, denoted by $J_s$ and $J_c$ respectively, as 
    \begin{align*}
        &J_s(f;\cD, \cD'):= \min \{k \in \mathbb{N} \cup \{0\}\;|\; \gamma_{f^{(k)}} + (-1)^{k+1}\, \bar\gamma_{f^{(k)}} \neq 0\},\quad\text{and}\\
        &J_c(f;\cD, \cD'):= \min \{k \in \mathbb{N} \cup \{0\}\;|\; \gamma_{f^{(k)}} + (-1)^k\, \bar\gamma_{f^{(k)}} \neq 0\}.
    \end{align*}
    We will drop the argument $(f;\cD, \cD')$ when it is clear from context.
\end{definition}
The signal index $J_s$ is the first $k$ such that an odd power of $\vu$ has non-zero coefficient in the decomposition of Proposition~\ref{prop:blockwigner} and the constant index $J_c$ is the the first $k$ such that an even power of $\vu$ has non-zero coefficient.

Community detection is straightforward when $J_s = 0$. In this case, even if the two communities initially have weights with the same mean $\gamma = \bar\gamma$, the difference in the means of the transformed matrix, $\gamma_{f^{(0)}}$ and $\bar\gamma_{f^{(0)}}$, allows a spectral algorithm to recover the community structure. This illustrates how transformation can aid signal recovery. Therefore, we now assume that $J_s \neq 0$. Similar to Section~\ref{sec:signed-signal-recovery}, we study signal recovery in two cases: $J_s < J_c$ and $J_s > J_c$.

In the case where $J_s > J_c$, the top eigenvalue corresponds to the eigenvector $\mathbf{1}$ and the second top eigenvector of $\bY$ is aligned to $\boldsymbol{\vu}$, if $|\gamma - \bar\gamma|$ is large enough. This is shown in the following theorem. Here, for simplicity we assume that we have a balanced stochastic block model with each community having the same number of nodes (i.e., $\beta = 1/2$), which enables us to derive closed-form expressions for the critical values of the phase transition. The same technique can be used for the general $\beta$ which gives a qualitatively similar result (with phase transitions happening at different critical values).
\begin{theorem}
    \label{thm:Js>Jc}
    In the stochastic block model defined above, let ${{\gamma - \bar\gamma}} = 2 c_\lambda n^{\alpha-1/2}$ with $\alpha \in [0,1/2), c_\lambda \in \R$.  Assume that the conditions of Proposition~\ref{prop:blockwigner} hold. Let $(\gamma_2, \vu_2)$ be the second top eigenvalue and eigenvector pair of $\frac{1}{\sqrt{n}}f(\bA)$, and define
    \vspace{-0.2cm}
    \begin{align*}
        \kappa = \frac{c_\lambda^{J_s}(\gamma_{f^{(J_s)} + (-1)^{J_s+1}} \bar\gamma_{f^{(J_s)}})}{2 J_s!}.
    \end{align*}
        \vspace{-0.2cm}
    \begin{enumerate}
        \item [(a)] If $\alpha < \frac{J_s - 1}{2 J_s}$, we have  $\gamma_2 \overset{\mathrm{P}}{\to} \sqrt{2(\sigma^2+\bar\sigma^2)}$ and $\frac{|\langle \vu_2, \vu\rangle|}{\Vert\vu_2\Vert_2 \Vert\vu\Vert_2} \overset{\mathrm{P}}{\to} 0$.
        \item [(b)] If $\alpha = \frac{J_s - 1}{2 J_s}$, we have
        \begin{align*}
            \hspace{-1.3cm}\gamma_2 \overset{\mathrm{P}}{\to}  \begin{cases}
                \sqrt{2(\sigma^2+\bar\sigma^2)} & \kappa < \sqrt{\frac{\sigma^2+\bar\sigma}{2}}\\
               \kappa +  {\frac{\sigma^2 + \bar\sigma^2}{2\kappa}}  & \kappa > \sqrt{\frac{\sigma^2+\bar\sigma}{2}}
            \end{cases}
            \quad\text{ and }\quad \frac{|\langle \vu_2, \vu\rangle|}{\Vert\vu_2\Vert_2 \Vert\vu\Vert_2} \overset{\mathrm{P}}{\to} \begin{cases}
                0 & \kappa < \sqrt{\frac{\sigma^2+\bar\sigma}{2}}\\
                 \sqrt{1 - \frac{2\kappa^2}{\sigma^2 + \bar\sigma^2}} & \kappa > \sqrt{\frac{\sigma^2+\bar\sigma}{2}}
            \end{cases}
        \end{align*}
        \item [(c)] If $\alpha > \frac{J_s - 1}{2J_s}$, we have  $\gamma_2 = \omega_\sP(1)$ and $\frac{|\langle \vu_2, \vu\rangle|}{\Vert\vu_2\Vert_2 \Vert\vu\Vert_2} \overset{\mathrm{P}}{\to} 1$.
    \end{enumerate} 
\end{theorem}
When $J_s < J_c$, if the signal strength is large enough, the top eigenvector will be aligned to $\boldsymbol{\zeta}$.
\begin{theorem}
    \label{thm:Js<Jc}
    Let conditions of Theorem~\ref{thm:IO>IE} hold, but with  $J_s < J_c$.
    Let $(\gamma_1, \vu_1)$ be the top eigenvalue and eigenvector pair of $\bY$. Then, Part (a), (b), and (c) of Theorem~\ref{thm:IO>IE} will hold for $\gamma_1, \vu_1$.
\end{theorem}

The proofs of these two theorems can be found in Section~\ref{sec:thm:Js>Jc} and Section~\ref{sec:thm:Js<Jc}. These theorems show that in order for it to be possible to recover the community membership using spectral methods, the signal strength $\sqrt{n}(\gamma - \bar\gamma)$ should grow $\Theta(n^{\frac{J_s - 1}{2J_s}})$ that is $\omega(1)$ if $J_s > 1$.

\section{Experiments}
\label{sec:experiments}

\paragraph{5.1. Signed Signal Recovery.}
\label{sec:experiment-signed}
Let $\bY$ follow \eqref{eq:main-form} with $\bW \sim \mathrm{Wig}_n(\normal(0,1))$, $\vx = \boldsymbol{\zeta}/\sqrt{n}$ with $\boldsymbol{\zeta}$ defined in Condition~\ref{cond:rademacher}, and $f(x) = (x^2 - 1) + (x^3 - 3x)$. For this setting, we have $I_o = 3$, and $I_e = 2$. Let $(\gamma_1, \vu_1)$ and $(\gamma_2, \vu_2)$ be the first and second top eigenvalue and eigenvector  of $\bY$.
\begin{figure}[h!]
     \centering
     \includegraphics[width=\textwidth]{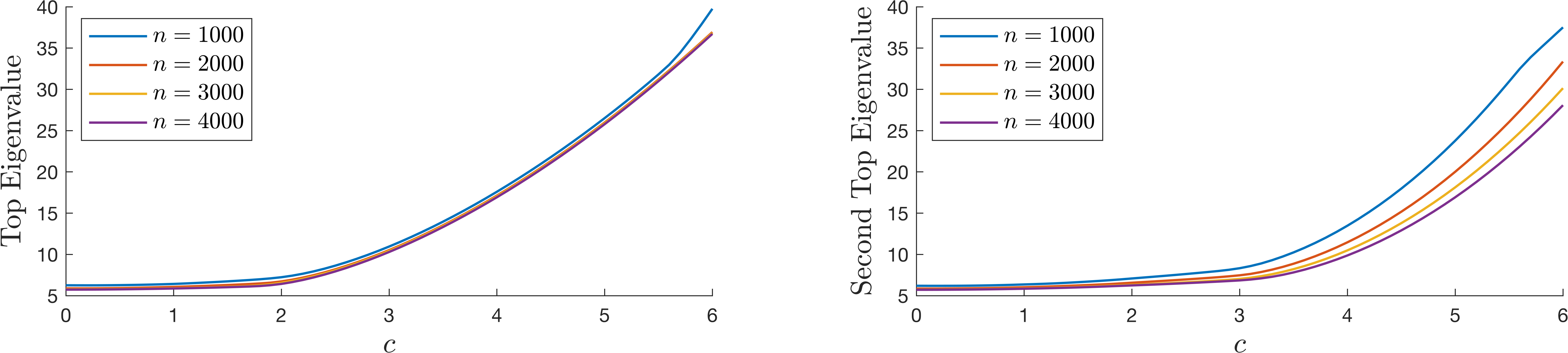}
     \caption{Position of the top and the second top eigenvalue of $\bY$ for the setting in Section~\ref{sec:experiment-signed}.1 with signal strength $\lambda = cn^{1/4}$ as a function of $c$ for different values of $n$.}
     \label{fig:signed-position}
 \end{figure}

\begin{figure}[h!]
     \centering
     \hspace{-0.5cm}
     \includegraphics[width=\textwidth]{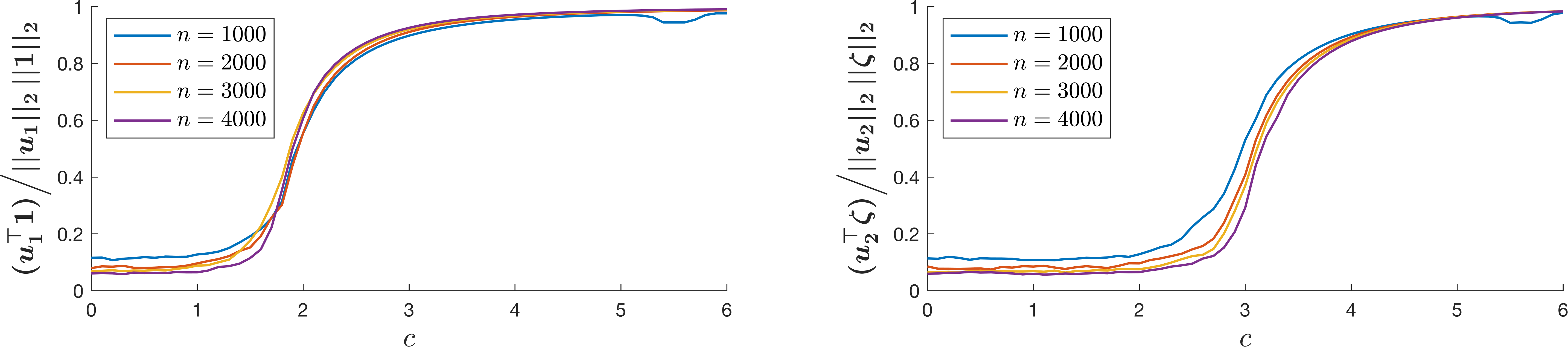}
     \caption{The correlation of the top eigenvector $\vu_1$ of $\bY$, and $\mathbf{1}$ and the correlation of the second eigenvector $\vu_2$ of $\bY$ and $\boldsymbol{\zeta}$ in the setting of Section~\ref{sec:experiment-signed}.1 with signal strength $\lambda = cn^{1/4}$ as a function of $c$ for different values of $n$.}
     \label{fig:signed-corr}
 \end{figure}

In Figure~\ref{fig:signed-position}, we set \(\lambda = c n^{1/4}\) and plot \(\gamma_1\) and \(\gamma_2\) as functions of \(c\) for \(n \in \{1000, 2000, 3000, 4000\}\). Figure~\ref{fig:signed-position} (Right) shows that as \(n\) varies, the curve for the position of \(\gamma_2\) shifts, whereas in Figure~\ref{fig:signed-position} (Left) the curve for the position of \(\gamma_1\) remains nearly unchanged. Also, in Figure~\ref{fig:signed-corr}, we plot the correlation of \(\vu_1\) with \(\mathbf{1}/\sqrt{n}\) and the correlation of \(\vu_2\) with \(\boldsymbol{\zeta}/\sqrt{n}\) as functions of \(c\). In Figure~\ref{fig:signed-corr} (Left), for all values of \(n\), the phase transition in the correlation of \(\vu_1\) with \(\mathbf{1}/\sqrt{n}\) occurs at the same threshold. However, in Figure~\ref{fig:signed-corr} (Right), the threshold for the phase transition of the correlation of \(\vu_2\) with \(\boldsymbol{\zeta}/\sqrt{n}\) shifts when $n$ varies.

This indicates that \(\lambda = \Theta(n^{1/4})\) is the appropriate scaling for the emergence of the first outlying eigenvalue but not for the second outlier. These findings align with the theoretical results from~Section~\ref{sec:signed-signal-recovery}, which suggest that for the emergence of the second top outlier, the signal strength should scale as \(\Theta(n^{1/3})\) when in the case where \(I_o = 3\). Additionally, comparing the left and the right figure, note that the phase transition for the second top eigenvector happens at a larger signal strength, demonstrating that signed signal recovery is more challenging than unsigned recovery in this problem.

\begin{figure}[h!]
    \centering
    \hspace{-0.5cm}\includegraphics[width=\textwidth]{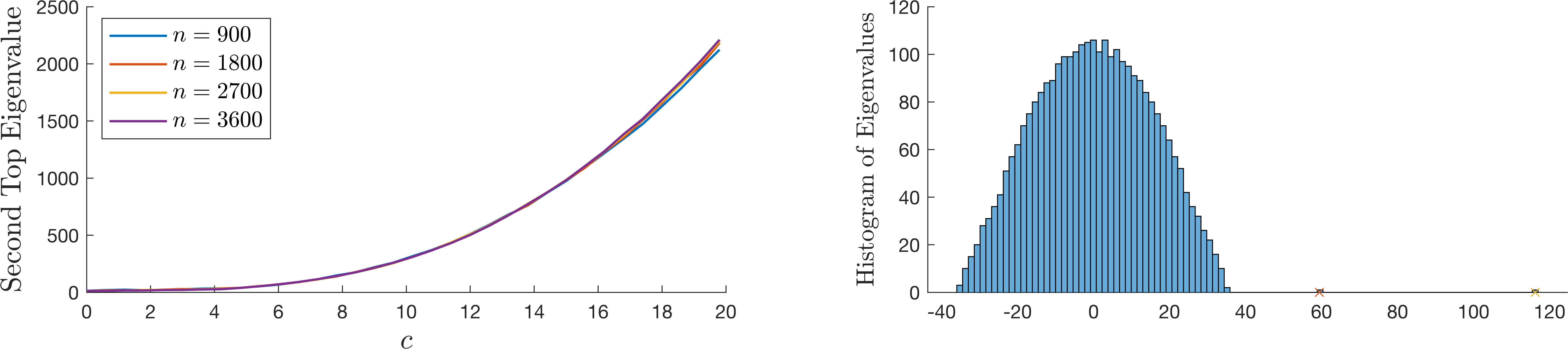}
    \caption{(Left) the position of the second top eigenvlaue of $\bY$ in the setting of Section~5.2. The curve is unchanged for different values of $n$. (Right) Histogram of the eigenvalues of $\bY$ in the setting of Section~5.2. The eigenvalues consist of a bulk of eigenvalues that stick together resulting from the noise component, and  two outliers corresponding to signal terms. }
    \label{fig:sbm}
\end{figure}

\paragraph{5.2. Transformed Stochastic Block Models.}
In the stochastic block model of Section~\ref{sec:pre-transformed-block}, let  $\cD = \normal(\Delta/2, 1)$, $\bar\cD = \normal(-\Delta/2, 0.5)$ for $\Delta \in \R$, and $\beta = 1/3$. Consider the nonlinear function $f(x) = 2.25(x^2 - 1) + (x^3 - 3x) + (x^4 - 6x^2 +3)$. This choice of parameters ensure that we have $\gamma - \bar\gamma = \Delta$, $\gamma_{f^{(0)}} - \bar \gamma_{f^{(0)}} = 0$, $J_s = 3$ and $J_c = 2$. Following the scaling in Theorem~\ref{thm:Js>Jc}, we set $\Delta = 2c n^{-1/6}$  (i.e., $\alpha = 1/3$). In Figure~\ref{fig:sbm} (Left), we plot the position of the top and eigenvalue of $\bY$ as functions of $c$. We see that  as \(n\) varies, the curve for the position of \(\gamma_1\) is unchanged, confirming that $\alpha = 1/3$ is the correct scale of signal strength for the emergence of the second top outlier in the spectrum of $\bY$ when $J_s = 3$. In Figure~\ref{fig:sbm} (Right), we plot the histogram of the eigenvalues of $\bY$. This histogram consists of a bulk of eigenvalues that stick together from the noise term in Theorem~\ref{thm:Js>Jc}, and two outlier signal components.

\section{Conclusion}
In this paper, we proposed a novel signal-plus-noise decomposition for a nonlinear spiked random matrix model in which a nonlinear function is applied element-wise to a rank-one perturbation of a random matrix. Through the study of the problems of signed signal recovery and community detection in transformed stochastic block models, we  showed that this novel decomposition enables the study of phase transitions in problems that could not be handled by prior works, and reveals interesting phenomena in these models. Nonlinear spiked random matrix model appears naturally in the study of feature learning in shallow neural networks (see e.g., \cite{ba2022high,moniri_atheory2023,dandi2023learning,cui2024asymptotics,wang2024nonlinear, dandi2024benefits}, etc.). We leave the potential applications to problems in deep learning theory as future work.

\section*{Acknowledgement}
Behrad Moniri and Hamed Hassani are supported by The Institute for Learning-enabled Optimization at Scale (TILOS), under award number NSF-CCF-2112665.

\bibliography{full_references}
\bibliographystyle{alpha}

\newpage
\appendix

\section{Proofs}

\subsection{Proof of Theorem~\ref{thm:sig+noise}}
\label{sec:thm:approximation-proof}
\begin{proof}
Consider the matrix model defined in \eqref{eq:main-form} and let $\lambda = c_\lambda n^\alpha$ with $\alpha \in [\frac{\ell-1}{2\ell},\frac{\ell}{2\ell+2})$. To analyze this model, we first consider a Taylor approximation of $\bY =\frac{1}{\sqrt{n}} f(\bW + \lambda \sqrt{n}\vx \vx^\top)$ around $\bW$; i.e., we write
\begin{align*}
     \bY = \frac{1}{\sqrt{n}}f(\bW) + \sum_{k = 1} ^{\ell} \frac{\lambda^k n^{(k-1)/2}}{k!}\left(f^{(k)}(\bW)\right) \circ \left(\vx^{\circ k}\vx^{\circ k\top}\right) + \cE,
\end{align*}
in which $f^{(k)}$ is the $k$-th derivative of $f$ and is applied element-wise.

We  show that the matrix $\cE$ will have a negligible operator norm.
Using the Taylor remainder theorem, we have
\begin{align*}
    \cE = \frac{\lambda^{\ell+1}n^{\frac{\ell}{2}}}{(\ell+1)!}\left(f^{(\ell+1)}(\bW+\bDelta) \right)\circ \left(\vx^{\circ (\ell+1)}\vx^{\circ (\ell+1)\top}\right),
\end{align*}
where $|\bDelta| \leq  \sqrt{n} \,|\lambda\vx\vx^\top|$ element-wise. The reminder matrix $\cE$ can be written as
\begin{align*}
    \cE = \frac{\lambda^{\ell+1}n^{\frac{\ell}{2}}}{(\ell+1)!}\Bigg[\text{diag}\left(\vx^{\circ (\ell+1)}\right)\left(f^{(\ell+1)}(\bW+\bDelta) \right) \text{diag}\left(\vx^{\circ (\ell+1)}\right)\Bigg].
\end{align*}
Hence, the operator norm of the matrix $\cE$ can be bounded as
\begin{align*}
    \Vert \cE \Vert_{\rm op} &\leq \frac{|\lambda|^{\ell+1}n^{\frac{\ell}{2}}}{(\ell+1)!}M_\vx^{2\ell+2} \Vert f^{(\ell+1)}(\bW+\bDelta) \Vert_{\rm op},
\end{align*}
where $M_\vx = \Vert\vx\Vert_\infty$.
Also, we can upper bound the operator norm of $f^{(\ell+1)}(\bW+\bDelta)$ using the Frobenius norm, to get
\begin{align*}
    \Vert f^{(\ell+1)}(\bW+\bDelta) \Vert_{\rm op}^2 &\leq \Vert f^{(\ell+1)}(\bW+\bDelta) \Vert_{\rm Fr}^2 \\
    &= \sum_{i = 1}^{n}\sum_{j = 1}^{n}\left( f^{(\ell+1)}(W_{ij}+\Delta_{ij}) \right)^2\\
    &\leq \sum_{i = 1}^{n}\sum_{j = 1}^{n}\left( C_f + (W_{ij}+\Delta_{ij})^{m_f} \right)^2\\
    &\leq \sum_{i = 1}^{n}\sum_{j = 1}^{n}\left( 2C_f^2 + 4 W_{ij}^{2m_f}+ 4\Delta_{ij}^{2m_f} \right),
\end{align*}
where in the second inequality, we have used Condition~\ref{cond:poly-bound}. Also, 
\begin{align*}
    \Delta_{ij}^{2m_f} \leq \lambda^{2m_f}n ^{m_f}M_\vx^{4m_f} \leq c_\lambda^{2m_f} n^{(2 \alpha-1) m_f} \kappa^{4m_f}(n) = o(1),
\end{align*}
where the last equality holds as long as $\kappa(n) = o(n^{\frac{1-2\alpha}{4}})$.     
With this, the expected value of $\Vert \cE \Vert_{\rm op}^2$ can be upper bounded as 
\begin{align*}
    \E \Vert \cE \Vert_{\rm op}^2 \leq \frac{|\lambda|^{2\ell+2}n^{\ell}}{((\ell+1)!)^2}M_\vx^{4\ell+4} \sum_{i = 1}^{n}\sum_{j = 1}^{n}\left( 2C_f^2 + 4 \E W_{ij}^{2m_f}+ o(1) \right).
\end{align*}
Using the fact that $\E W_{ij}^{2m_f}$ is bounded using Condition~\ref{cond:strength}, we have
\begin{align*}
    \E \Vert \cE \Vert_{\rm op} = O_\sP\left(\lambda^{\ell+1}n^{\frac{\ell+2}{2}} M_\vx^{2\ell+2}\right) = O_\sP\left(n^{\left(\alpha - \frac{1}{2}\right) \ell+\alpha} \kappa^{2\ell+2}(n)\right).
\end{align*}
Hence, $\E\Vert\cE\Vert_{\rm op} = o(1)$ because $\alpha < \frac{\ell}{2\ell+2}$.

Then, using Markov's inequality, we have $\Vert\cE\Vert_{\rm op} = o(1)$ with probability $1 - o(1)$. Putting everything together, we have
\begin{align*}
    \bY = \frac{1}{\sqrt{n}  }f(\bW) + \sum_{k = 1} ^{\ell} \frac{\lambda^k n^{(k-1)/2}}{k!}\left( f^{(k)}(\bW)\right) \circ \left(\vx^{\circ k}\vx^{\circ k\top}\right) + \cE,
\end{align*}
with $\Vert\cE\Vert_{\rm op} = o_\sP(1)$. Next, we show that the matrix
\begin{align*}
    \mathbf{\Gamma} := \sum_{k = 1} ^{\ell} \frac{\lambda^k n^{(k-1)/2}}{k!}\left(  f^{(k)}(\bW) - \E f^{(k)}(\bW)\right) \circ \left(\vx^{\circ k}\vx^{\circ k\top}\right)
\end{align*}
has $\Vert\mathbf{\Gamma}\Vert_{\rm op} = o_\sP(1)$. To show this, we write
\begin{align*}
    \Vert\mathbf{\Gamma}\Vert_{\rm op} &\leq \sum_{k = 1} ^{\ell} \frac{|\lambda|^k n^{(k-1)/2}}{k!} \left\Vert\left(  f^{(k)}(\bW) - \E f^{(k)}(\bW)\right) \circ \left(\vx^{\circ k}\vx^{\circ k\top}\right)\right\Vert_{\rm op}\\
    &\leq \sum_{k = 1} ^{\ell} \frac{|\lambda|^k n^{(k-1)/2}M_\vx^{2k}}{k!} \left\Vert f^{(k)}(\bW) - \E f^{(k)}(\bW)\right\Vert_{\rm op}.
\end{align*} 
Also, the matrix $f^{(k)}(\bW) - \E f^{(k)}(\bW)$ is a mean zero Wigner-type random matrix and its operator norm is known to be $O_\sP(\sqrt{n})$ (see e.g., \cite[Theorem 2.1]{erdHos2019bounds}). A similar order-wise argument as the one used above gives $\Vert\mathbf{\Gamma}\Vert_{\rm op} = o_\sP(1)$. Hence, putting everything together, we get
\begin{align*}
    \bY = \frac{1}{\sqrt{n}}f(\bW) + \sum_{k = 1} ^{\ell} \frac{\lambda^k n^{(k-1)/2}}{k!}\left( \E f^{(k)}(\bW)\right) \circ \left(\vx^{\circ k}\vx^{\circ k\top}\right) + \tilde\cE,
\end{align*}
with $\Vert\tilde\cE\Vert_{\rm op} = o_\sP(1)$.
\end{proof}
\subsection{Proof of Proposition~\ref{prop:wigner}}
\label{sec:thm:wigner}
When $\bW$ is a Wigner matrix, for any $k \in \mathbb{N}$, we have
\begin{align*}
    \E f^{(k)}(\bW) = \mu_{f^{(k)}}\mathbf{1}\mathbf{1}^\top.
\end{align*}
Under the assumptions of this proposition, Theorem~\ref{thm:sig+noise} holds with
\begin{align*}
    \bH_k = \frac{\lambda^k n^{(k-1)/2} \mu_{f^{(k)}}}{k!} \vx^{\circ k} \vx^{\circ k \top}.
\end{align*}
Also, denoting $\mathbf{\Lambda} = \frac{1}{\sqrt{n}} \left[f(\bW) - \E f(\bW)\right]$, we have
\begin{align*}
    \frac{1}{\sqrt{n}}f(\bW) &= \mathbf{\Lambda} + \frac{\mu_{f^{(0)}}}{\sqrt{n}}\mathbf{1}\mathbf{1}^\top.
\end{align*}
Here, $\mathbf{\Lambda}$ mean-zero, and it is a symmetric matrix with i.i.d. entries (up to symmetric constraints), thus it is a mean-zero Wigner matrix. These, alongside with Theorem~\ref{thm:sig+noise} conclude the proof.

\subsection{Proof of Theorem~\ref{thm:IO>IE}}
\label{sec:thm:IO>IE}
Under the assumptions of this theorem, using \eqref{eq:Wigner-Rad-decomposition} we have
\begin{align}
    \label{eq:signed-rec-first-decompose}
    \bY = \mathbf{\Lambda} +    \left(\sum_{k = 1, k \text{ odd}} ^{\ell} \frac{c_\lambda^k n^{(\alpha-1/2)k}\mu_{f^{(k)}}}{k! \sqrt{n}}\right)\boldsymbol{\zeta}\boldsymbol{\zeta}^{\top} + \left(\sum_{k = 0, k \text{ even}} ^{\ell} \frac{c_\lambda^k n^{(\alpha-1/2)k}\mu_{f^{(k)}}}{ k!\sqrt{n}}\right)\mathbf{1}\mathbf{1}^{\top} + \bDelta.
\end{align}
We will prove each case where $\alpha < \frac{I_o-1}{2I_o}$, $\alpha = \frac{I_o-1}{2I_o}$, and $\alpha > \frac{I_o-1}{2I_o}$ separately.

\paragraph{Case 1.} If $\alpha < \frac{I_o-1}{2I_o}$, based on \eqref{eq:signed-rec-first-decompose}, we have
\begin{align*}
    \bY = \mathbf{\Lambda} +  \frac{\kappa_1}{n} \mathbf{1}\mathbf{1}^\top
    + \bDelta,
\end{align*}
with $\Vert\bDelta\Vert_{\rm op} = o_\sP(1)$, and
\begin{align*}
    \kappa_1 = 
        \frac{c_\lambda^{I_e} n^{(\alpha-1/2){I_e}}\mu_{f^{({I_e})}}}{ {I_e}!\sqrt{n}}.
\end{align*}
In this case, the phase transition for the top eigenvalue of the spiked Wigner models in Section~\ref{sec:linear_wigner} readily proves that $\gamma_2 \overset{\mathrm{P}}{\to} 2\sigma_f$.  Using the Sherman–Morrison formula, we can write
\begin{align*}
    \frac{1}{n}\,\boldsymbol{\zeta}^\top (z\bI - \mathbf{\Lambda} -  \kappa_1 \mathbf{1}\mathbf{1}^\top/n)^{-1}\;\boldsymbol{\zeta}  &= \frac{1}{n}\,\boldsymbol{\zeta}^\top (z\bI - \mathbf{\Lambda})^{-1}\boldsymbol{\zeta} - \frac{\frac{\kappa_1}{n^2} (\boldsymbol{\zeta}^\top (z\bI - \mathbf{\Lambda})^{-1} \mathbf{1})^2}{1 + \frac{\kappa_1}{n} \mathbf{1}^\top (z\bI - \mathbf{\Lambda})^{-1}\mathbf{1}}\\
    &=  \frac{1}{n}\,\boldsymbol{\zeta}^\top (z\bI - \mathbf{\Lambda})^{-1}\boldsymbol{\zeta} + o_\sP(1),
\end{align*}
where in the last line we have used the fact that $\frac{1}{n}\E\,{\boldsymbol{\zeta}^\top\bG \mathbf{1}} = 0.$
Thus, using the Hanson-Wright inequality, we have 
\begin{align*}
    \frac{1}{n}{\boldsymbol{\zeta}^\top\bG \mathbf{1}}  = O_\sP(1/\sqrt{n}).
\end{align*}

On the other hand, defining $(\gamma_i, \vu_i)$ as the top $i$-th eigenvalue and eigenvector pair of $\bY$, we have
\begin{align*}
    \frac{1}{n}\,\boldsymbol{\zeta}^\top (z\bI - \mathbf{\Lambda} -  \kappa_1 \mathbf{1}\mathbf{1}^\top/n)^{-1} \,\boldsymbol{\zeta}= \frac{1}{n}\,\sum_{i = 1}^{n} (z - \gamma_i)^{-1} (\boldsymbol{\zeta}^\top \vu_i)^2.
\end{align*}
Putting these together and noting that $\gamma_2 \overset{\mathrm{P}}{\to} 2\sigma_f$, we can write
\begin{align*}
    \frac{1}{n}\,(\boldsymbol{\zeta}^\top \vu_2)^2 &= \lim_{z \to 2\sigma_f} \frac{z-2\sigma_f}{n}\,\boldsymbol{\zeta}^\top(z\bI - \mathbf{\Lambda} -  \kappa_1 \mathbf{1}\mathbf{1}^\top/n)^{-1}\,\boldsymbol{\zeta}\\
    &= \lim_{z \to 2\sigma_f} \frac{2\sigma_f - z}{2\sigma_f^2} \left(z - \sqrt{z^2 - 4\sigma_f^2}\right) = 0,
\end{align*}
proving the theorem for this case. Here, we have used the classic random matrix theory result (see e.g., \cite{pizzo2013finite}, etc.) that
\begin{align*}
    \frac{1}{n}\, \boldsymbol{\zeta}^\top (z\bI - \mathbf{\Lambda} )^{-1} \, \boldsymbol{\zeta} \overset{P}{\to} m_{\rm sc}(z),
\end{align*}
where $m_{\rm sc}(z)$ is the Stieltjes transform of the semicircle law and is given by
\begin{align*}
    m_{\rm sc}(z) = \frac{-1}{2\sigma_f^2}\left(z - \sqrt{z^2 - 4\sigma_f^2}\right).
\end{align*}

\paragraph{Case 2.} If $\alpha = \frac{I_o-1}{2I_o}$ (and noting that $I_e < I_o$), we have
\begin{align*}
\bY = \mathbf{\Lambda} + 
    \frac{\kappa_1}{n} \mathbf{1}\mathbf{1}^{\top} + \frac{\kappa_2}{n}\boldsymbol{\zeta}\boldsymbol{\zeta}+ \bDelta,    
\end{align*}
in which $\kappa_1$ and $\kappa_2$ are given by
\begin{align*}
    &\kappa_1 = \frac{c_\lambda^{I_e} \mu_{f^{({I_e})}}}{I_e! } \left({n^{\frac{1}{2} - \frac{I_e}{2I_o}}} + o_\sP\left({n^{\frac{1}{2} - \frac{I_e}{2I_o}}} \right)\right) = \Theta_\sP\left({n^{\frac{1}{2} - \frac{I_e}{2I_o}}} \right) = \omega_\sP(1),\\
    &\kappa_2 = \frac{c_\lambda^{I_o} \mu_{f^{({I_o})}}}{I_o!} = \Theta_\sP(1).
\end{align*}

The eigenvalues of $\bY$ are the solutions to $\det(\bY - z\bI) = 0$. We write
\begin{align*}
    \det(z\bI - \bY) &= \det\left( z\bI - \mathbf{\Lambda} - \kappa_1 \mathbf{1}\mathbf{1}^\top/n - \kappa_2 \boldsymbol{\zeta}\boldsymbol{\zeta}^{\top}/n  \right)\\
    &=\det\left((z\bI - \mathbf{\Lambda} - \kappa_1 \mathbf{1}\mathbf{1}^\top/n)^{-1}  \right)\det\left(\bI + \kappa_2(z\bI - \mathbf{\Lambda} - \kappa_1 \mathbf{1}\mathbf{1}^\top/n)^{-1}\boldsymbol{\zeta}\boldsymbol{\zeta}^{\top}/n  \right) = 0,
\end{align*}
where we have assumed that $z$ is not an eigenvalue of $\mathbf{\Lambda} + \kappa_1 \mathbf{1}\mathbf{1}^\top/n$. Thus, we have
\begin{align*}
    \det\left(\bI + \kappa_2( z\bI -\mathbf{\Lambda} - \kappa_1 \mathbf{1}\mathbf{1}^\top/n)^{-1}\boldsymbol{\zeta} \boldsymbol{\zeta}^\top/n \right) =0
\end{align*}
which implies that
\begin{align*}
    \label{eq:fixed-eq}
    \frac{1}{n}\,\boldsymbol{\zeta}^\top(z \bI - \mathbf{\Lambda} - \kappa_1 \mathbf{1}\mathbf{1}^\top/n)^{-1} \boldsymbol{\zeta} = - \frac{1}{\kappa_2}.
\end{align*}
Next, using the Sherman–Morrison formula, we write
\begin{align}
    (z\bI - \mathbf{\Lambda} - \kappa_1 \mathbf{1}\mathbf{1}^\top/n)^{-1} = \bG - \frac{\kappa_1\,\bG\,\mathbf{1}\mathbf{1}^\top\bG/n}{1 + \kappa_1\, \mathbf{1}^\top \bG \mathbf{1}/n}
\end{align}
where $\bG = (z\bI - \mathbf{\Lambda})^{-1}$. Hence,
\begin{align*}
    \frac{1}{n}\,\boldsymbol{\zeta}^\top(  z\bI 
 -\mathbf{\Lambda} - \kappa_1 \mathbf{1}\mathbf{1}^\top/n)^{-1} \boldsymbol{\zeta} = \frac{1}{n}\,\boldsymbol{\zeta}^\top\bG\boldsymbol{\zeta}- \frac{\kappa_1\,\left(\boldsymbol{\zeta}^\top\bG\,\mathbf{1}\right)^2}{n^2(1 + \kappa_1\, \mathbf{1}^\top \bG \mathbf{1}/n)}.
\end{align*}
Note that $\frac{1}{n}{\boldsymbol{\zeta}^\top\bG\boldsymbol{\zeta}}$ and $\frac{1}{n}\mathbf{1}^\top \bG\, \mathbf{1}$ are both $\Theta_\sP(1)$, and $\frac{1}{n}\E\,{\boldsymbol{\zeta}^\top\bG \mathbf{1}} = 0.$
Thus, using the Hanson-Wright inequality, we have 
\begin{align*}
    \frac{1}{n}{\boldsymbol{\zeta}^\top\bG \mathbf{1}}  = O_\sP(1/\sqrt{n}).
\end{align*}
Putting everything together, we get
\begin{align*}
    \frac{\kappa_1\,\left(\boldsymbol{\zeta}^\top\bG\,\mathbf{1}\right)^2}{n(1 + \kappa_1\, \mathbf{1}^\top \bG \mathbf{1}/n)}\overset{\mathrm{P}}{\to} 0.
\end{align*}
As a result, from \eqref{eq:fixed-eq}, for any $z$ that is an eigenvalue of $\bY$ and is not an eigenvalue of $\boldsymbol{\Lambda} + \kappa_1 \mathbf{1}\mathbf{1}^\top/n$ we have 
\begin{align}
    \label{eq:cite-later-for-the-sbm}
    \frac{1}{n}\boldsymbol{\zeta}^\top(z\bI 
 -\mathbf{\Lambda})^{-1} \boldsymbol{\zeta} = - \frac{1}{\kappa_2} + o_\sP(1).
\end{align}
Using the Hanson-Wright inequality, we can prove that
\begin{align*}
    \frac{1}{n}\,\boldsymbol{\zeta}^\top(z\bI - \mathbf{\Lambda})^{-1} \boldsymbol{\zeta} =\frac{1}{n}\,\E \,{\boldsymbol{\zeta}^\top(z\bI - \mathbf{\Lambda})^{-1} \boldsymbol{\zeta}} + o_\sP(1)   \overset{\mathrm{P}}{\to} \frac{1}{n}\mathrm{tr}\left((z\bI -\mathbf{\Lambda})^{-1}\right).
\end{align*}
It is known that when $\mathbf{\Lambda}$ is a Wigner matrix, the trace of the matrix $(z\bI - \mathbf{\Lambda})^{-1}$ will converge in probability to the Stieltjes transform of the semi-circle law (see e.g., \cite{pizzo2013finite}, etc.); i.e., 
\begin{align*}
 \frac{1}{n} \mathrm{tr}\left((z\bI - \mathbf{\Lambda} )^{-1}\right) \overset{\mathrm{P}}{\to} -\frac{1}{2\sigma_f^2} \left(z - \sqrt{z^2 - 4\sigma_f^2}\right).   
\end{align*}
Thus, the eigenvalue $z$ should asymptotically satisfy
\begin{align*}
    \frac{1}{2\sigma_f^2} \left(z - \sqrt{z^2 - 4\sigma_f^2}\right) =  \frac{1}{\kappa_2}
\end{align*}
and solving this equation assuming $\kappa_2 > \sigma_f$, gives
\begin{align}
    \label{eq:eigen_in_proof}
    z = \kappa_2 + \frac{\sigma_f^2}{\kappa_2},
\end{align}
and if $\kappa_2 < \sigma_f$, it has no solution and the second top eigenvalue is the edge of the limiting empirical eigenvalue distribution of $\mathbf{\Lambda}$. Thus we have
\begin{align}
    \label{eq:limit_gamma_2}
    \gamma_2 \overset{P}{\to} \begin{cases}
        2\sigma_f & \kappa_2 < \sigma_f\\
        \kappa_2 + \frac{\sigma_f^2}{\kappa_2} & \kappa_2 > \sigma_f.
    \end{cases}
\end{align}

Next, we study the alignment between the second top eigenvector $\vu_2$ and $\boldsymbol{\zeta}$.  Similar to the first case, defining $\bG = (z \bI - \boldsymbol{\Lambda})^{-1}$, we can use the Woodbury formula to get
\begin{align*}
    \left(z\bI - \mathbf{\Lambda} - \frac{\kappa_1}{n}\mathbf{1}\mathbf{1}^\top - \frac{\kappa_2}{n}\boldsymbol{\zeta}\boldsymbol{\zeta}^\top \right)^{-1} = \bG - \bG\,\frac{\left(\frac{1}{\kappa_1}+ \frac{\mathbf{1}^\top\bG\mathbf{1}}{n}\right) \frac{\boldsymbol{\zeta}\boldsymbol{\zeta}^\top}{n} + \left(\frac{1}{\kappa_2}+ \frac{\boldsymbol{\zeta}^\top\bG\boldsymbol{\zeta}}{n}\right) \frac{\boldsymbol{1}\boldsymbol{1}^\top}{n} + \bDelta}{\left(\frac{1}{\kappa_2} + \frac{\boldsymbol{\zeta}^\top \bG\boldsymbol{\zeta}}{n} \right)\left(\frac{1}{\kappa_1} + \frac{\boldsymbol{1}^\top \bG\boldsymbol{1}}{n} \right) + o_\sP(1)} \,\bG,
\end{align*}
where $\gamma_i, \vu_i$ is the $i$-th top eigenvalue and eigenvector pair of $\bY$ and  $\Vert\bDelta\Vert_{\rm op} = o_\sP(1)$.  On the other hand, we have
\begin{align*}
    \frac{1}{n}\,\boldsymbol{\zeta}^\top \left(z\bI - \mathbf{\Lambda} - \frac{\kappa_1}{n}\mathbf{1}\mathbf{1}^\top - \frac{\kappa_2}{n}\boldsymbol{\zeta}\boldsymbol{\zeta}^\top \right)^{-1} \,\boldsymbol{\zeta}= \frac{1}{n}\,\sum_{i = 1}^{n} (z - \gamma_i)^{-1} (\boldsymbol{\zeta}^\top \vu_i)^2.
\end{align*}
Thus, noting that $\gamma$ will converge in probability to a deterministic value $\bar\gamma_2$ (with its exact expression given in \eqref{eq:limit_gamma_2}) we can write
\begin{align*}
    \frac{1}{n}(\boldsymbol{\zeta}^\top \vu_2)^2 &= \lim_{z \to \bar\gamma_2} \frac{z - \bar\gamma_2}{n}\,\boldsymbol{\zeta}^\top \left(z\bI - \mathbf{\Lambda} - \frac{\kappa_1}{n}\mathbf{1}\mathbf{1}^\top - \frac{\kappa_2}{n}\boldsymbol{\zeta}\boldsymbol{\zeta}^\top \right)^{-1}  \,\boldsymbol{\zeta}\\[0.2cm]
    &= \lim_{z \to \bar\gamma_2} \frac{z - \bar\gamma_2}{n} \left(\boldsymbol{\zeta}^\top  \bG \,\boldsymbol{\zeta} - \frac{ {(\boldsymbol{\zeta}^\top\bG\,\boldsymbol{\zeta})^2}/{n^2} + o_\sP(1)}{{1}/{\kappa_2} + {\boldsymbol{\zeta}^\top \bG\,\boldsymbol{\zeta}}/{n}}\right)\\[0.2cm]
    &= \lim_{z \to \bar\gamma_2} {(z - \bar\gamma_2)}\left(m_{\rm sc}(z) - \frac{ m^2_{\rm sc}(z)}{{1}/{\kappa_2} + m_{\rm sc}(z)}\right).
\end{align*}
Here, we have again used the fact that $\frac{1}{n}{\boldsymbol{\zeta}^\top\bG \mathbf{1}}  = O_\sP(1/\sqrt{n})$ and that $\frac{1}{n}{\boldsymbol{\zeta}^\top\bG \mathbf{1}}  = O_\sP(1/\sqrt{n})$ and that $\frac{1}{n}\boldsymbol{\zeta}^\top\bG \boldsymbol{\zeta} \overset{P}{\to} m_{\rm sc}(z)$. 

The first term is equal does not contribute to the limit and is equal to zero. From \eqref{eq:limit_gamma_2}, when $\kappa_2 < \sigma_f$, we have $\bar\gamma_2 = 2\sigma_f$; thus, $(\boldsymbol{\zeta}^\top \vu_2)^2/n \overset{P}{\to} 0$. However, when $\kappa_2 > \sigma_f$, we have that $\bar\gamma_2 $ satisfies $m_{\rm sc}(\bar\gamma_2) = -1/\kappa_2$ and using the L'Hôpital's rule, we get
\begin{align*}
    \frac{1}{n}(\boldsymbol{\zeta}^\top \vu_2)^2 \overset{P}{\to} \frac{m^2_{\rm sc}(\bar\gamma_2)}{m'_{\rm sc}(\bar\gamma_2)}.
\end{align*}
Plugging in the formula for $m_{\rm sc}$ and also the explicit expression for $\bar\gamma_2$ in \eqref{eq:limit_gamma_2}, we show
\begin{align*}
    \frac{1}{n}(\boldsymbol{\zeta}^\top \vu_2)^2 \overset{P}{\to}  1- \frac{\sigma_f^2}{\kappa_2^2},
\end{align*}
completing the proof in this case.

\paragraph{Case 3.}  If $\alpha > \frac{I_o-1}{2I_o}$, we have
\begin{align*}
\bY = \mathbf{\Lambda} + 
    \frac{\kappa_1}{n} \mathbf{1}\mathbf{1}^{\top} + \frac{\kappa_2}{n}\boldsymbol{\zeta}\boldsymbol{\zeta}+ \bDelta,    
\end{align*}
with $\kappa_2= \omega_\sP(1)$ and $\kappa_1 = \omega_\sP(\kappa_2)$. The proof in this case is completed by applying the Davis-Kahan theorem \cite[Theorem 4.5.5]{vershynin2018high} to the second top eigenvector of $\bY$, and by noting $\kappa_2$ is asymptotically separated from $\kappa_1$ and $\Vert\mathbf{\Lambda}\Vert_{\rm op}$.

\subsection{Proof of Theorem~\ref{thm:IO<IE}}
\label{sec:thm:IO<IE}
Under the assumptions of this theorem, using \eqref{eq:Wigner-Rad-decomposition} we have
\begin{align*}
    \bY = \mathbf{\Lambda} +    \kappa\,\boldsymbol{\zeta}\boldsymbol{\zeta}^{\top} + \bDelta.
\end{align*}
where $\mathbf{\Lambda}$ is a Wigner matrix, $\Vert\bDelta\Vert_{\rm op} = o_\sP(1)$, and
\begin{align*}
    \kappa = \frac{c_\lambda^{I_o} n^{(\alpha-1/2){I_o}}\mu_{f^{({I_o})}}}{{I_o}! \sqrt{n}} + o_\sP\left(\frac{c_\lambda^{I_o} n^{(\alpha-1/2){I_o}}\mu_{f^{({I_o})}}}{{I_o}! \sqrt{n}}\right).
\end{align*}
This shows that in this setting, the matrix $\bY$ can be approximated by a spiked Wigner model. The proof is then completed by using Proposition~\ref{prop:wigner} for the phase transition in (linear) spiked Wigner models.

\subsection{Proof of Proposition \ref{prop:blockwigner}}
\label{sec:prop:blockwigner}
The observed matrix in this problem can be written as
$f(\bA) = f(\bar\bA + \frac{\gamma - \bar\gamma}{2} \vu \vu^\top)$ where $\bar\bA = \bA - \E \bA$ which means that $f(\bA)$ is a nonlinear spiked random matrix model defined in \eqref{eq:main-form} with $\bW = \bar\bA$, $\vx = \vu/\sqrt{n}$, and $\lambda = \frac{{\gamma - \bar\gamma}}{2} \sqrt{n}$. As a result, based on Theorem~\ref{thm:sig+noise}, we have
\begin{align*}
    \frac{1}{\sqrt{n}}f(\bA) =  \frac{1}{\sqrt{n}}f(\bar\bA) + \sum_{k = 1} ^{\ell} \frac{c_\lambda^k n^{(\alpha-1/2)k}}{k!\sqrt{n}}\left( \E f^{(k)}(\bar\bA)\right) \circ \left(\vu^{\circ k}\vu^{\circ k\top}\right) +\mathbf{\Delta},
\end{align*}
where $\Vert\mathbf{\Delta}\Vert_{\rm op} = o_\sP(1)$. Note that the deterministic matrix $\E f^{(k)}(\bar\bA)$ is block-structured and can be written as
\begin{align*}
    \E f^{(k)}(\bar\bA) = \begin{bmatrix}
        \gamma_{f^{(k)}}\,\mathbf{1}_{\beta n \times \beta n} & \bar\gamma_{f^{(k)}}\,\mathbf{1}_{\beta n \times (1-\beta)n}\\
        \bar\gamma_{f^{(k)}}\,\mathbf{1}_{(1-\beta)n \times \beta n)} & \gamma_{f^{(k)}}\,\mathbf{1}_{(1-\beta)n \times (1-\beta)n} 
    \end{bmatrix} = \frac{ \gamma_{f^{(k)}} + \bar\gamma_{f^{(k)}}}{2} \mathbf{1}\mathbf{1}^\top + \frac{  \gamma_{f^{(k)}} - \bar\gamma_{f^{(k)}}}{2} \boldsymbol{u}\boldsymbol{u}^\top.
\end{align*}

Plugging this into the expansion, we arrive at  $\Vert\bar\bY - \frac{1}{\sqrt{n}}f(\bA)\Vert_{\rm op} = o_\sP(1)$ for
\begin{align*}
     \bar\bY = \frac{1}{\sqrt{n}}f(\bar \bA) + \sum_{k = 1} ^{\ell} \frac{c_\lambda^k n^{(\alpha-1/2)k}}{k! \sqrt{n}} \left(\frac{ \gamma_{f^{(k)}} + \bar\gamma_{f^{(k)}}}{2}\vu^{\circ k}\vu^{\circ k\top}  + \frac{ \gamma_{f^{(k)}} - \bar\gamma_{f^{(k)}}}{2} \vu^{\circ (k+1)}\vu^{\circ (k+1)\top}\right).
\end{align*}
The first term has a mean equal to 
\begin{align*}
    \E f(\bar\bA) = \frac{ \gamma_{f^{(0)}} + \bar\gamma_{f^{(0)}}}{2} \mathbf{1}\mathbf{1}^\top + \frac{  \gamma_{f^{(0)}} - \bar\gamma_{f^{(0)}}}{2} \boldsymbol{u}\boldsymbol{u}^\top. 
\end{align*}
Adding as subtracting this term, we arrive at 
\begin{align*}
     \bar\bY = \mathbf{\Gamma} + \sum_{k = 0} ^{\ell} \frac{c_\lambda^k n^{(\alpha-1/2)k}}{k! \sqrt{n}} \left(\frac{ \gamma_{f^{(k)}} + \bar\gamma_{f^{(k)}}}{2}\vu^{\circ k}\vu^{\circ k\top}  + \frac{ \gamma_{f^{(k)}} - \bar\gamma_{f^{(k)}}}{2} \vu^{\circ (k+1)}\vu^{\circ (k+1)\top}\right),
\end{align*}
which concludes the proof of the theorem.

\subsection{Proof of Proposition \ref{prop:wigner-type-block}}
\label{sec:prop:wigner-type-block}
Note that the matrix $\mathbf{\Gamma} = \frac{1}{\sqrt{n}} (f(\bar\bA) - \E f(\bar\bA))$ is mean zero and it has independent entries (up to symmetric constraints). Thus, it satisfies Definition~\ref{def:wigner-type} and is a Wigner-type matrix. The matrix of variances for this matrix ensemble is 
\begin{align*}
    [\bS]_{ij} = \E |[\mathbf{\Gamma}]_{ij}|^2 = \frac{1}{n}\begin{bmatrix}
        \sigma^2 \mathbf{1}_{\beta n \times \beta n}&\bar\sigma^2 \mathbf{1}_{\beta n \times (1-\beta) n}\\
        \bar\sigma^2 \mathbf{1}_{\beta n \times \beta n}&\sigma^2 \mathbf{1}_{\beta n \times (1-\beta) n} \in \R^{n\times n}.
    \end{bmatrix}
\end{align*}
The matrix $\mathbf{\Gamma}$ satisfies the conditions of Theorem~\ref{theorem:local_law}. The QVE in \eqref{eq:QVE} is symmetric for all $i \in \{1, \dots, \beta n\}$ and also for all $i \in \{\beta n+1, \dots, n\}$. Denoting $\bar m_1 = m_1 = \dots = m_{\beta n}$ and $\bar m_2 = m_{\beta n +1 } = \dots = m_{n}$, \eqref{eq:QVE} can be simplified as
\begin{align*}
    &-\frac{1}{\bar{m}_1(z)} = z + \beta \sigma^2 \bar m_1(z) + (1-\beta)\bar\sigma^2\bar m_2(z), \\
    &-\frac{1}{\bar{m}_2(z)} = z + \beta \bar\sigma^2 \bar m_2(z) + (1-\beta)\sigma^2\bar m_1(z).
\end{align*}
Having this, Theorem~\ref{theorem:local_law} concludes the proof.

\subsection{Proof of Theorem \ref{thm:Js>Jc}}
\label{sec:thm:Js>Jc}
Based on Proposition~\ref{prop:blockwigner}, we have
\begin{align*}
     \bar\bY = \mathbf{\Gamma} + \sum_{k = 0} ^{\ell} \frac{c_\lambda^k n^{(\alpha-1/2)k}}{k! \sqrt{n}} \left(\frac{ \gamma_{f^{(k)}} + \bar\gamma_{f^{(k)}}}{2}\vu^{\circ k}\vu^{\circ k\top}  + \frac{ \gamma_{f^{(k)}} - \bar\gamma_{f^{(k)}}}{2} \vu^{\circ (k+1)}\vu^{\circ (k+1)\top}\right).
\end{align*}
We will prove each case where $\alpha < \frac{J_s-1}{2J_s}$, $\alpha = \frac{J_s-1}{2J_s}$, and $\alpha > \frac{J_s-1}{2J_s}$ separately.

\paragraph{Case 1.} When $\alpha < \frac{J_s-1}{2J_s}$, 
\begin{align*}
    \bar\bY = \mathbf{\Gamma} + \frac{\kappa_c}{n} \mathbf{1}\mathbf{1}^\top + \mathbf{\Delta},
\end{align*}
in which $\Vert\bDelta\Vert_{\rm op} = o_\sP(1)$, and 
\begin{align*}
    \kappa_c = 
         \frac{c_\lambda^{J_c} n^{(\alpha-1/2)J_c + 1/2}( \gamma_{f^{(J_c)}} +  (-1)^{J_c}\, \bar\gamma_{f^{(J_c)}})}{2 J_c!},
\end{align*}
with $\Vert\bDelta\Vert_{\rm op} = o_\sP(1)$.  In this case, and assuming that $\beta = 1/2$, we readily have $\gamma_2 \overset{\mathrm{P}}{\to} \bar\gamma_2 := \sqrt{2(\sigma^2 + \bar\sigma^2)}$ which is the right edge of the limiting empirical eigenvalue distribution of $\mathbf{\Gamma}$. Following case 1 in the proof of Theorem~\ref{thm:IO>IE},  we can write
\begin{align*}
    \frac{1}{n}\,(\boldsymbol{u}^\top \vu_2)^2 &= \lim_{z \to \bar\gamma_2 } \frac{z-\bar\gamma_2 }{n}\;\boldsymbol{u}^\top \left(z\bI - \mathbf{\Gamma} -  \frac{\kappa_c}{n} \mathbf{1}\mathbf{1}^\top\right)^{-1}\;\boldsymbol{u} .
\end{align*}
Using the Sherman–Morrison formula, we can write
\begin{align*}
    \frac{1}{n}\,\boldsymbol{u}^\top \left(z\bI - \mathbf{\Gamma} -  \frac{\kappa_c}{n} \mathbf{1}\mathbf{1}^\top\right)^{-1}\boldsymbol{u}  &= \frac{1}{n}\,\boldsymbol{u}^\top (z\bI - \mathbf{\Gamma})^{-1}\boldsymbol{u} - \frac{\frac{\kappa_c}{n^2} (\boldsymbol{u}^\top (z\bI - \mathbf{\Gamma})^{-1} \mathbf{1})^2}{1 + \frac{\kappa_c}{n} \mathbf{1}^\top (z\bI - \mathbf{\Gamma})^{-1}\mathbf{1}}
\end{align*}
Note that according to the second part of Theorem~\ref{theorem:local_law}, we have $\boldsymbol{u}^\top (z\bI - \mathbf{\Gamma})^{-1}\boldsymbol{u}/n \overset{P}{\to} \bar m(z)$ and $\boldsymbol{1}^\top (z\bI - \mathbf{\Gamma})^{-1}\boldsymbol{1}/n \overset{P}{\to} \bar m(z)$ 
 where $\bar m(z)$ is given in Proposition~\ref{prop:wigner-type-block}. Also, using the same theorem, $\vu^\top (z\bI - \boldsymbol{\Gamma})^{-1}\mathbf{1}/n \overset{P}{\to} 0$. Using these, we can take the limit $z \to \bar\gamma_2 $ to get
\begin{align*}
    \frac{1}{n}\,(\boldsymbol{u}^\top \vu_2)^2\overset{P}{\to} 0,
\end{align*}
completing the proof in this case.

\paragraph{Case 2.} If $\alpha = \frac{J_s-1}{2J_s}$, we have
\begin{align*}
    \bar\bY = \mathbf{\Gamma} + \frac{\kappa_s}{n} \vu\vu^\top + \frac{\kappa_c}{n} \mathbf{1}\mathbf{1}^\top + \mathbf{\Delta},
\end{align*}
where $\Vert\bDelta\Vert_{\rm op} = o_\sP(1)$, and
\begin{align*}
    &\kappa_s = \frac{c_\lambda^{J_s} n^{(\alpha-1/2)J_s + 1/2}( \gamma_{f^{(J_s)}} +  (-1)^{J_s+1}\, \bar\gamma_{f^{(J_s)}})}{2 J_s!} = \Theta(1),\\
    &\kappa_c = \frac{c_\lambda^{J_c} n^{(\alpha-1/2)J_c + 1/2}( \gamma_{f^{(J_c)}} +  (-1)^{J_c}\, \bar\gamma_{f^{(J_c)}})}{2 J_c!} = \omega(1).
\end{align*}


Similar to the proof in Section~\ref{sec:thm:IO>IE} for case 2, note that the eigenvalues of $\bY$ are the solutions to $\det(\bY - z\bI) = 0$. Following the same line of reasoning used to prove \eqref{eq:cite-later-for-the-sbm}, and noting that here we also have $\vu^\top (z\bI - \boldsymbol{\Gamma})^{-1}\mathbf{1}/n \overset{P}{\to} 0$, we get
\begin{align*}
    \frac{1}{n}\,\vu^\top(\mathbf{\Gamma} - z\bI)^{-1} \vu = - \frac{1}{\kappa_s} + o_\sP(1).
\end{align*}
Using this and the second part of Theorem~\ref{theorem:local_law}, we have 
\begin{align*}
    \frac{1}{n}\,\vu^\top(\mathbf{\Gamma} - z\bI)^{-1} \vu \overset{\mathrm{P}}{\to} \bar m(z) 
\end{align*}
Hence, the outlying eigenvalue is the solution to the nonlinear equation
\begin{align}
    \label{eq:nonline-eq}
    \bar m(z) = -\frac{1}{\kappa_s}.
\end{align}
Again, note that when $\beta = 1/2$, the solution to the quadratic vector equation of Proposition~\ref{prop:blockwigner} is given by
\begin{align*}
    \bar m(z) = -\frac{1}{2\sigma_f^2} \left(z - \sqrt{z^2 - 4\sigma_f^2}\right),
\end{align*}
with $\sigma_f = \sqrt{\frac{\sigma^2 + \bar\sigma^2}{2}}$. Using this to solve \eqref{eq:nonline-eq}, we get
\begin{align}
    \gamma_2 \overset{P}{\to} \begin{cases}
        \sqrt{2(\sigma^2 + \bar\sigma^2)} & \kappa_s < \sqrt{\frac{\sigma^2 + \bar\sigma^2}{2}}\\
        \kappa_s +  {\frac{\sigma^2 + \bar\sigma^2}{2\kappa_s}} & \kappa_s > \sqrt{\frac{\sigma^2 + \bar\sigma^2}{2}}.
    \end{cases}
\end{align}


To analyze the eigenvector alignment, we follow the same line of reasoning as in case 2 of the proof of Theorem~\ref{thm:IO>IE}: 
\begin{align*}
    \frac{1}{n}(\vu^\top \vu_2)^2 &= \lim_{z \to \bar\gamma_2} \frac{z - \bar\gamma_2}{n}\,\vu^\top \left(z\bI - \mathbf{\Lambda} - \frac{\kappa_\mathbf{1}}{n}\mathbf{1}\mathbf{1}^\top - \frac{\kappa_s}{n}\vu\vu^\top \right)^{-1}  \,\vu\\[0.2cm]
    &= \lim_{z \to \bar\gamma_2} \frac{z - \bar\gamma_2}{n} \left(\vu^\top  \bG \,\vu - \frac{ {(\vu^\top\bG\,\vu)^2}/{n^2} + o_\sP(1)}{{1}/{\kappa_s} + {\vu^\top \bG\,\vu}/{n}}\right)\\[0.2cm]
    &= \lim_{z \to \bar\gamma_2} {(z - \bar\gamma_2)}\left(\bar m(z) - \frac{ \bar m^2(z)}{{1}/{\kappa_s} + \bar m(z)}\right).
\end{align*}
Here, we have used the fact that $\vu^\top (z\bI - \boldsymbol{\Gamma})^{-1}\mathbf{1}/n \overset{P}{\to} 0$. When $\kappa_s < \sqrt{\frac{\sigma^2 + \bar\sigma^2}{2}}$, we have $\bar\gamma_2 = \sqrt{2(\sigma^2 + \bar\sigma^2)}$; thus, $(\boldsymbol{v}^\top \vu_2)^2/n \overset{P}{\to} 0$. However, when $\kappa_s > \sqrt{\frac{\sigma^2 + \bar\sigma^2}{2}}$, we have that $\bar\gamma_2 $ satisfies $\bar m(\bar\gamma_2) = -1/\kappa_s$ and using the L'Hôpital's rule, we get
\begin{align*}
    \frac{1}{n}(\boldsymbol{v}^\top \vu_2)^2 \overset{P}{\to} \frac{\bar m^2(\bar\gamma_2)}{\bar m'(\bar\gamma_2)} = 1 - \frac{2\kappa_s^2}{\sigma^2 + \bar\sigma^2}> 0,
\end{align*}
which completes the proof of the theorem.

\paragraph{Case 3.}  If $\alpha > \frac{J_s-1}{2J_s}$, the proof is completed by applying the Davis-Kahan theorem \cite[Theorem 4.5.5]{vershynin2018high} to the second top eigenvector of $\bY$, and by noting $\kappa_s$ is asymptotically separated from $\kappa_c$ and the bulk eigenvalues $\Vert\mathbf{\Lambda}\Vert_{\rm op}$. 

\subsection{Proof of Theorem \ref{thm:Js<Jc}}
\label{sec:thm:Js<Jc}
Based on Proposition~\ref{prop:blockwigner}, we have
\begin{align*}
     \bar\bY = \mathbf{\Gamma} + \sum_{k = 0} ^{\ell} \frac{c_\lambda^k n^{(\alpha-1/2)k}}{k! \sqrt{n}} \left(\frac{ \gamma_{f^{(k)}} + \bar\gamma_{f^{(k)}}}{2}\vu^{\circ k}\vu^{\circ k\top}  + \frac{ \gamma_{f^{(k)}} - \bar\gamma_{f^{(k)}}}{2} \vu^{\circ (k+1)}\vu^{\circ (k+1)\top}\right).
\end{align*}
We will prove each case where $\alpha < \frac{J_s-1}{2J_s}$, $\alpha = \frac{J_s-1}{2J_s}$, and $\alpha > \frac{J_s-1}{2J_s}$ separately.

\paragraph{Case 1.} When $\alpha < \frac{J_s-1}{2J_s}$, we can write
\begin{align*}
    \bar\bY  =\mathbf{\Gamma} + \bDelta,
\end{align*}
with $\Vert\bDelta\Vert_{\rm op} = o_\sP(1)$. This readily proves that $\gamma_1 \overset{P}{\to} \sqrt{2(\sigma^2 + \bar\sigma^2)}$ which is the right edge of the limiting empirical eigenvalue distribution of $\mathbf{\Gamma}$. The corresponding eigenvector can be shown to have vanishing correlation with $\vu$ by using the eigenvector delocalization result from \cite[Corollary 1.14]{ajanki2017universality}.

\paragraph{Case 2.} When $\alpha = \frac{J_s-1}{2J_s}$, we can write
\begin{align*}
    \bar\bY = \mathbf{\Gamma} + \frac{\kappa_s}{n} \vu\vu^\top + \mathbf{\Delta},
\end{align*}
in which $\Vert\bDelta\Vert_{\rm op} = o_\sP(1)$, and 
\begin{align*}
    \kappa_s = \frac{c_\lambda^{J_s} ( \gamma_{f^{(J_s)}} +  (-1)^{J_s}\, \bar\gamma_{f^{(J_s)}})}{2 J_s!} = \Theta(1).
\end{align*}

Similar to the proof in Section~\ref{sec:thm:IO>IE} for case 2, note that the eigenvalues of $\bar\bY$ are the solutions to $\det(\bar\bY - z\bI) = 0$. Following the same line of reasoning used to prove \eqref{eq:cite-later-for-the-sbm}, we get
\begin{align*}
    \frac{1}{n}\,\boldsymbol{u}^\top(\mathbf{\Gamma} - z\bI)^{-1} \boldsymbol{u} = - \frac{1}{\kappa_s} + o_\sP(1).
\end{align*}
Using this proposition and the second part of Theorem~\ref{theorem:local_law}, we have 
\begin{align*}
    \frac{1}{n}\,\boldsymbol{u}^\top(\mathbf{\Gamma} - z\bI)^{-1} \boldsymbol{u} \overset{\mathrm{P}}{\to} \bar m(z).
\end{align*}
Hence, the outlying eigenvalue is the solution to the nonlinear equation
\begin{align}
    \label{eq:nonline-eq2}
    \bar m(z) = -\frac{1}{\kappa_s}.
\end{align}
Solving this equation, noting that when $\beta = 1/2$, the solution to the quadratic vector equation of Proposition~\ref{prop:blockwigner} is given by
\begin{align*}
    \bar m(z) = -\frac{1}{2\sigma_f^2} \left(z - \sqrt{z^2 - 4\sigma_f^2}\right),
\end{align*}
with $\sigma_f = \sqrt{\frac{\sigma^2 + \bar\sigma^2}{2}}$, we get
\begin{align}
    \label{eq:limit_gamma_3}
    \gamma_1 \overset{P}{\to} \begin{cases}
        \sqrt{2(\sigma^2 + \bar\sigma^2)} & \kappa_s < \sqrt{\frac{\sigma^2 + \bar\sigma^2}{2}}\\
        \kappa_s +  {\frac{\sigma^2 + \bar\sigma^2}{2\kappa_s}} & \kappa_s > \sqrt{\frac{\sigma^2 + \bar\sigma^2}{2}}.
    \end{cases}
\end{align}

To analyze the alignment of the outlying eigenvector and $\vu$, we follow the same line of reasoning as case 2 in the proof of Theorem~\ref{thm:IO>IE} and write
\begin{align*}
    \frac{1}{n}\,(\boldsymbol{u}^\top \vu_1)^2 &= \lim_{z \to \bar\gamma_1} \frac{z-\bar\gamma_1}{n}\;\boldsymbol{u}^\top \left(z\bI - \mathbf{\Gamma} -  \frac{\kappa_s}{n} \vu\vu^\top\right)^{-1}\;\boldsymbol{u} .
\end{align*}
Using the Sherman–Morrison formula, we can write
\begin{align*}
    \frac{1}{n}\,\boldsymbol{u}^\top \left(z\bI - \mathbf{\Gamma} -  \frac{\kappa_s}{n} \vu\vu^\top\right)^{-1}\;\boldsymbol{u}  &= \frac{1}{n}\,\boldsymbol{u}^\top (z\bI - \mathbf{\Gamma})^{-1}\boldsymbol{u} - \frac{\frac{\kappa_s}{n^2} (\boldsymbol{u}^\top (z\bI - \mathbf{\Gamma})^{-1} \vu)^2}{1 + \frac{\kappa_s}{n} \vu^\top (z\bI - \mathbf{\Gamma})^{-1}\vu}
\end{align*}
Putting these together, and using the fact that $\vu^\top(z\bI - \mathbf{\Gamma})^{-1}\vu \overset{P}{\to} \bar m(z)$, we get 
\begin{align*}
    \frac{1}{n}(\boldsymbol{u}^\top \vu_1)^2 = \lim_{z \to \bar\gamma_2} {(z - \bar\gamma_1)}\left(\bar m(z) - \frac{ \bar m^2(z)}{{1}/{\kappa_s} + \bar m(z)}\right).
\end{align*}
The first term is equal does not contribute to the limit and is equal to zero. Let $\bar\gamma_1 = \lim_{n \to \infty}\gamma_1$. From \eqref{eq:limit_gamma_3}, when $\kappa_s < \sqrt{\frac{\sigma^2 + \bar\sigma^2}{2}}$, we have $\bar\gamma_1 = \sqrt{2(\sigma^2 + \bar\sigma^2)}$; thus, $(\boldsymbol{\zeta}^\top \vu_1)^2/n \overset{P}{\to} 0$. However, when $\kappa_s > \sqrt{\frac{\sigma^2 + \bar\sigma^2}{2}}$, we have that $\bar\gamma_1 $ satisfies $\bar m(\bar\gamma_1) = -1/\kappa_s$ and using the L'Hôpital's rule, we get
\begin{align*}
    \frac{1}{n}(\boldsymbol{\zeta}^\top \vu_1)^2 \overset{P}{\to} 1 - \frac{2\kappa_s^2}{\sigma^2 + \bar\sigma^2}> 0,
\end{align*}
completing the proof.

\paragraph{Case 3.}  If $\alpha > \frac{J_s-1}{2J_s}$, the proof is completed by applying the Davis-Kahan theorem \cite[Theorem 4.5.5]{vershynin2018high} to the second top eigenvector of $\bY$, and by noting the rank one components and the bulk are asymptotically separated.

\section{Wigner-Type Random Matrices}
\label{sec:wigner-type}
First, we recall the definition of Wigner-type matrices from the random matrix theory literature (see e.g., \cite{ajanki2017universality,erdos2019matrix}, etc.).

\begin{definition}[Wigner-type matrix]
\label{def:wigner-type}
Let $\bH \in \R^{n\times n}$ be a real symmetric random matrix. The matrix $\bH$ is said to be Wigner-type if its entries $H_{ij}$ are mean zero and independent (up to symmetric constraints).
\end{definition}
Next, for a Wigner-type ensemble, the matrix of variances $\bS = [S_{ij}] \in \R^{n \times n}$ is defined as
\begin{align*}
S_{i j}=\mathbb{E}\left|H_{i j}\right|^2.
\end{align*}
For every such matrix $\bS$, the \textit{quadratic vector equation} (QVE), 
\begin{align}
    \label{eq:QVE}
    -\frac{1}{m_i(z)} = z + \sum_{j = 1}^{n} S_{ij} m_j(z), \quad \text{for all}\quad i \in [n], \quad z \in \mathbb{C}_+,
\end{align}
for $\mathbf{m} = (m_1, . . . , m_n) : \mathbb{C}_+ \to \mathbb{C}_+^{n}$ on the complex upper half   plane has a unique solution \cite{ajanki2017universality}. The case where the matrix $\bS$ is doubly stochastic is referred to as the generalized Wigner matrix, and $m_i(z) = m_{\rm sc}(z); \forall i \in [n]$ where $m_{\rm sc}(z)$ is the Stieltjes transform of the semi-circle law, is the solution to \eqref{eq:QVE} (see, e.g., \cite{erdos2019matrix}, etc.).

Next, we state the local law for Wigner-type matrices from \cite{ajanki2017universality}, i.e. that for large $n$, under some conditions on $\bS$, the resolvent $\bG(z) = (z\mathbf{I} - \bH)^{-1}$, with spectral parameter $z = \tau + i\eta \in \mathbb{C}_+$, is close to $\mathrm{diag}(\vm(z))$, as long as $\eta \gg \frac{1}{n}$. To do this, first, we recall the notion of stochastic domination.
\begin{definition}[Stochastic domination] Suppose $n_0:\R_+^2 \rightarrow \mathbb{N}$. For two sequences of random variables $\varphi^{(n)}$ and $\psi^{(n)}$, of non-negative random variables we say that $\varphi^{(n)}$ is stochastically dominated by $\psi^{(n)}$ if for all $\varepsilon>0$ and $D>0$,
$$
\mathbb{P}\left(\varphi^{(n)}>n^{\varepsilon} \psi^{(n)}\right) \leq n^{-D}, \quad n \geq n_0(\varepsilon, D) .
$$
In this case we write $\varphi^{(n)} \prec \psi^{(n)}$.

\end{definition}

We are now ready to state the local law for the Wigner-type random matrices. Let $\mathbf{m} = (m_1, . . . , m_n)$ be the unique solution to the QVE and define the density
\begin{align}
    \label{eq:rho-from-stieltjes}
    \rho(\tau) := \lim_{\eta \searrow 0} \frac{1}{\pi n } \sum_{i = 1}^{n} \mathbf{Imag} \;m_i(\tau + i \eta).
\end{align}
The properties of $\rho$ are studied in \cite[Theorem 4.1 and Corollary 1.3]{ajanki2017universality}. The following result is proved in \cite[Theorems 1.7 and 1.13]{ajanki2017universality}.
\begin{theorem}[Local Law for Wigner-Type Ensembles]
    \label{theorem:local_law}
    Let $\bH \in \R^{n\times n}$ be a Wigner-type ensemble with entries that have bounded moments, and assume that for fixed $a, b \in \R$ and $L \in \mathbb{N}$, we have that for any $n$,
            \begin{align*}
            S_{ij} \leq \frac{a}{n}, \quad \text{and} \quad (\bS^L)_{ij} \geq \frac{b}{n}, \quad i,j \in [n].
        \end{align*}
    Fix some $\gamma>0$. Then, for any $z = \tau + i \eta\in \mathbb{C}_+$ with $\eta \geq n^{\gamma - 1}$, the resolvent matrix $\bG(z) = (\bH- z\mathbf{I})^{-1}$ satisfies
    \begin{align*}
        \max_{i,j\in [n]} \left|G_{ij}(z) - m_i(z) \delta_{ij}\right| \prec \frac{1 + \sqrt{\rho(z)}}{\sqrt{\eta n}} + \frac{1}{\eta n}.
    \end{align*}
    Also, for any two deterministic unit $\ell_2$-norm vectors $\vw, \vv \in \R^n$,   
    \begin{equation*}
    \left|\sum_{i, j=1}^n {w}_i G_{i j}(z) v_j-\sum_{i=1}^n m_i(z) {w}_i v_i\right| \prec \frac{1 + \sqrt{\rho(z)}}{\sqrt{\eta n}} + \frac{1}{\eta n}.
    \end{equation*}
\end{theorem}






\end{document}

%% file: defs.tex
\usepackage{amsmath,amsfonts,bm}

\def\eqref#1{equation~\ref{#1}}

\def\1{\bm{1}}

\def\vm{{\bm{m}}}

\def\vu{{\bm{u}}}
\def\vv{{\bm{v}}}
\def\vw{{\bm{w}}}
\def\vx{{\bm{x}}}

\DeclareMathAlphabet{\mathsfit}{\encodingdefault}{\sfdefault}{m}{sl}
\SetMathAlphabet{\mathsfit}{bold}{\encodingdefault}{\sfdefault}{bx}{n}

\def\sP{{\mathbb{P}}}

\newcommand{\E}{\mathbb{E}}

\newcommand{\R}{\mathbb{R}}

\newcommand{\bY}{\mathbf{Y}}
\newcommand{\bW}{\mathbf{W}}

\newcommand{\bB}{\mathbf{B}}
\newcommand{\bA}{\mathbf{A}}
\newcommand{\bS}{\mathbf{S}}

\newcommand{\bH}{\mathbf{H}}

\newcommand{\bG}{\mathbf{G}}

\newcommand{\cV}{\mathcal{V}}
\newcommand{\cE}{\mathcal{E}}
\newcommand{\cG}{\mathcal{G}}

\newcommand{\bDelta}{\mathbf{\Delta}}

\newcommand{\bI}{\mathbf{I}}
\newcommand{\normal}{{\sf N}}

\newcommand{\cD}{\mathcal{D}}

\newcommand*\colourcheck[1]{%
  \expandafter\newcommand\csname #1check\endcsname{\textcolor{#1}{\ding{52}}}%
}

\colourcheck{green}

\usepackage{amsthm}
\newtheorem{theorem}{Theorem}[section]
\newtheorem{proposition}[theorem]{Proposition}

\newtheorem{definition}[theorem]{Definition}

\newtheorem{condition}[theorem]{Condition}

\renewcommand{\epsilon}{\varepsilon}

\usepackage{amsthm, amsopn, bm, xcolor, amsmath}
\usepackage{multirow}
\usepackage{booktabs}
\usepackage{pifont}

\usepackage{microtype}
\usepackage{graphicx}
\usepackage{subfigure}
\usepackage{booktabs} 

\usepackage{amssymb}
\usepackage{mathtools}
\usepackage[utf8]{inputenc}

\def\hmath$#1${\texorpdfstring{{\rmfamily\textit{#1}}}{#1}}

\usepackage{setspace}